\documentclass[final,3p,times]{elsarticle}

\usepackage{amssymb,amsfonts,amsmath}
\usepackage{amsthm}
\usepackage{algorithmic}
\usepackage{algorithm}
\usepackage{array}
\usepackage{booktabs}
\usepackage{caption} 
\usepackage{empheq}
\usepackage{graphicx}
\usepackage{lipsum}
\usepackage{mathrsfs}
\usepackage{multirow}
\usepackage{textcomp}
\usepackage{stfloats}
\usepackage{subfigure}
\usepackage{url}
\usepackage{verbatim}
\usepackage[active]{srcltx}
\usepackage{tabularx}
\theoremstyle{definition}
\newtheorem{theorem}{Theorem}
\newtheorem{lemma}{Lemma}
\newtheorem{definition}{Definition}
\newtheorem{corollary}{Corollary}
\newtheorem{proposition}{Proposition}

\newtheorem{remark}{Remark}

\DeclareMathOperator{\sgn}{sgn}

\newcommand{\tp}{^\top}

\newcommand{\lng}{\langle}
\newcommand{\rng}{\rangle}
\newcommand{\lf}{\left}
\newcommand{\rg}{\right}

\newcommand{\R}{\mathbb R}

\newcommand{\B}{\mathbb B}

\newcommand{\bnr}{{\cal B}_v(n,r)}

\newcommand{\snr}{{\rm St}(n,r)}

\newcommand{\Rnr}{\mathbb{R}^{n\times r}}

\journal{Signal Processing}

\begin{document}

\begin{frontmatter}



\title{A relaxation method for binary orthogonal optimization problems with its applications}

\author[label1]{Lianghai Xiao}
\author[label2]{Yitian Qian \corref{cor2}}
\author[label3]{and Shaohua Pan}
\cortext[cor2]{Corresponding author.\\
   Email address:  xiaolh$@$jnu.edu.cn (L. Xiao), yitian.qian@polyu.edu.hk(Y. Qiao), shhpan$@$scut.edu.cn  (S. Pan). }
 
\affiliation[label1]{organization={Jinan University},
             city={Guangzhou},
             country={China}}
\affiliation[label2]{organization={Hong Kong Polytechnic University},
             city={Hong Kong},
             country={China}}
\affiliation[label3]{organization={South China University of Technology},
             city={Guangzhou},
             country={China}}

\begin{abstract}
This paper focuses on a class of binary orthogonal optimization problems frequently arising in semantic hashing. Consider that this class of problems may have an empty feasible set, rendering them not well-defined. We introduce an equivalent model involving a restricted Stiefel manifold and a matrix box set, and then investigate its penalty problems induced by the $\ell_1$-distance from the box set and its Moreau envelope. The two penalty problems are always well-defined. Moreover, they serve as the global exact penalties provided that the original feasible set is non-empty. Notably, the penalty problem induced by the Moreau envelope is a smooth optimization over an embedded submanifold with a favorable structure. We develop a retraction-based line-search Riemannian gradient method to address the penalty problem. Finally, the proposed method is applied to supervised and unsupervised hashing tasks and is compared with several popular methods on the MNIST and CIFAR-10 datasets. The numerical comparisons reveal that our algorithm is significantly superior to other solvers in terms of feasibility violation, and it is comparable even superior to others in terms of evaluation metrics related to the Hamming distance. 
\end{abstract}




\end{frontmatter}



\section{Introduction}
\label{sec1}
Let $\mathbb{R}^{n\times r}$ be the vector space consisting of all $n\times r$ real matrices, equipped with the trace inner product and its induced Frobenious norm $\|\cdot\|_F$. We are concerned with the following binary orthogonal optimization problem
\begin{align}\label{eq:sh}
	&\min_{B \in\mathbb{R}^{n\times r}}\widetilde{f}(B), \nonumber\\
	&\quad{\rm s.t.}\ \  B\in\mathcal{M}(n,r)\cap \mathbb{L}_{v},
\end{align}
where $\widetilde{f}\!:\Rnr\rightarrow\R$ is a continuously differentiable function,
$\mathcal{M}(n,r)\!:=\{B\in\{-1,1\}^{n\times r}~|~B^{\top}B=nI_r\}$ is the binary orthogonal set, and $\mathbb{L}_{v}\!:=\big\{B\in\mathbb{R}^{n\times r}\,|\,B^{\top}v=0\big\}$ is a subspace associated to a given vector $v\in\mathbb{R}^{n}$. Here, $I_r$ denotes an identity matrix of dimension $r$.

Problem \eqref{eq:sh} often appears in semantic hashing  \cite{salakhutdinov2007learning}, a crucial technique to provide effective solutions for information retrieval and text classification. Here we present the objective functions of representative models categorized as binary orthogonal optimization problems \eqref{eq:sh}; see Table \ref{tab:ref}.

\begin{table*}[ht]
\centering
\caption{Representative models of the binary orthogonal optimization}\label{tab:ref}
\footnotesize
\begin{tabularx}{\textwidth}{ll|llll}
\toprule
NO. & Model & Year & $v$ & Objective function & Method \\\midrule
1 & SH\cite{weiss2008spectral} & 2008 & $e$ & $\min_B{\rm tr}(B^\top LB)$ & Relaxation \\
2 & STH\cite{zhang2010self} & 2010 & $e$ & $\min_B{\rm tr}(B^\top LB)$ & Relaxation \\
3 & AGH\cite{liu2011hashing} & 2011 & $e$ & $\min_B{\rm tr}(B^\top LB)$ & Relaxation \\
4 & ITQ\cite{gong2012iterative} & 2012 & $0$ & $\min_{B,R}\|B-XWR\|_F^2$ & Alternating maximization \\
5 & KSH\cite{liu2012supervised} & 2012 & 0 & $\min_{B}\|\frac{1}{r}BB^\top-S\|_F^2$ & Relaxation \\
6 & MSH\cite{weiss2012multidimensional} & 2012 & $e$ & $\min_{B,R}\|W-BLB^\top\|_F^2$ & Relaxation \\
7 & DGH\cite{liu2014discrete} & 2014 & $e$ & $\min_B{\rm tr}(B^\top LB)$ & Relaxation,  alternating maximization \\
8 & CH\cite{liu2014collaborative} & 2014 & 0 & $\min_{B,R}\frac{1}{n}\|B-RW^\top\|_F^2+\frac{\lambda}{n}\|B^\top G - U\|_F^2$ & Relaxation \\
9 & RDSH\cite{yang2015robust} & 2015 & 0 & $\min_{B,R,Q}{\rm tr}(R^\top LR) + \lambda\|B - RQ\|_F^2$ & Alternating maximization \\
10 & BDNN\cite{do2016learning} & 2016 & $e$ & $\min_{B,R}\|Y-BR^\top\|_F^2 + \frac{\lambda}{2}\|R\|_F^2$ & Alternating maximization \\
11 & DPLM\cite{shen2016fast} & 2016 & $e$ & $\min_B{\rm tr}(B^\top LB)$ or $\min_{B,R}\frac{1}{2}\|Y-B^\top R\|_F^2 + \lambda\|R\|_F^2$ & Penalty, projected gradient \\
12 & SHSR\cite{hu2018discrete} & 2018 & $e$ & $\min_{B,R}\frac{1}{2}\|YR-B\|_F^2 + \lambda\|R\|_F^2$ & Relaxation, alternating maximization \\
13 & DCSH\cite{hoang2020unsupervised} & 2020 & $e$ & $\min_{B,\{R_m\}_M}\sum_{m=1}^M\left({\rm tr}(R_m^\top L R_m) - \lambda{\rm tr}(B^\top R_m)\right)$ & Penalty, binary optimization \\
14 & RSLH\cite{liu2020reinforced} & 2020 & 0 & $\min_{B, R, Q}\|Y - R^\top B\|_F^2 + \lambda\|B - Q^\top Y\|_F^2$ & Alternating maximization \\
15 & SCDH\cite{chen2020strongly} & 2020 & $e$ & $\min_{B}\|YS-BB^\top\|_F^2$ & Alternating maximization \\
16 & ABMO\cite{xiong2021generalized} & 2021 & $e$ &  $\min_B{\rm tr}(B^\top LB)$ or $\min_{B,R}\frac{1}{2}\|Y-B^\top R\|_F^2 + \lambda\|R\|_F^2$ & ADMM \\
17 & UDH\cite{jin2021unsupervised} & 2021 & $e$ & $\min_{B,R}{\rm tr}(B^\top L B) +\lambda\|B-YR\|_F^2$ & Relaxation \\
18 & SCLCH\cite{qin2022joint} & 2022 & $e$ & $\min_{B,\{R_m\}_M}\sum_{m=1}^{M}\|B^\top R_m -rL^\top L\|_F^2$ & Alternating maximization \\
19 & HHL\cite{sun2023hierarchical} & 2023 & $e$ & $\min_{B,R,Q}\|B-YRQ\|_F^2+\lambda\|RQ\|_{2,1}^2$ & Relaxation \\
\bottomrule
\end{tabularx}
\end{table*}

 A notable limitation of problem \eqref{eq:sh} is the potential emptiness of the set $\mathcal{M}(n,r)$, rendering it not well-defined (for instance, when $n$ is odd). Therefore, specific relation between $n$ and $r$ is required to  ensure that $\mathcal{M}(n,r)\ne\emptyset$. It has been observed that those matrices in $\mathcal{M}(n,r)$ are actually the ``Orthogonal Arrays'' with level 2. For the definition of orthogonal arrays and the sufficient condition of non-emptiness of $\mathcal{M}(n,r)$, interested readers can find more details in \cite{hedayat1999orthogonal}.

 As validated by Lemma \ref{equiv-model} later, the feasible set of \eqref{eq:sh}  is equivalent to $\{B\in[-E,E]~|~ B\tp B = nI_r,B\tp v = 0\}$, where $[-E,E]$ denotes the matrix box set with $E$ representing an $n\times r$ real matrix of all ones. Consequently,  the binary orthogonal optimization model \eqref{eq:sh} can be equivalently reformulated as the following problem
 \begin{align}\label{eq:ori}
 &\min_{X \in\mathbb{R}^{n\times r}} f(X):=\widetilde{f}(\sqrt{n}X), \nonumber\\
 &\quad{\rm s.t.}~X\in\bnr\cap\Lambda,
 \end{align}
 in the sense that they share the same global and local optimal solution sets, where $\bnr\!:={\rm St}(n,r)\cap\mathbb{L}_{v}$ with ${\rm St}(n,r):=\{X\in\!\Rnr\,|\,X\tp X = I_r\}$, and $\Lambda\!:=[-\frac{1}{\sqrt{n}}E,\frac{1}{\sqrt{n}}E]$. We call $\bnr$ the restricted Stiefel manifold, a submanifold of the Stiefel manifold ${\rm St}(n,r)$. Clearly, the feasible set of \eqref{eq:ori} may also be empty. Instead of solving \eqref{eq:sh} or \eqref{eq:ori} directly, we focus on the penalty problem of \eqref{eq:ori} induced by the $\ell_1$-norm distance function to the constraint set $\Lambda$:
 \begin{equation}\label{eq:penprob}
 \underset{X\in\bnr}{\min}\ F(X):=f(X)+\rho h_{\frac{1}{\sqrt{n}}}(X),
 \end{equation}
 where $\rho>0$ is a parameter, and $h_{\frac{1}{\sqrt{n}}}(X):={\rm dist}_1(X,\Lambda)$ is the $\ell_1$-norm distance function onto the box set $\Lambda$. Problem \eqref{eq:penprob} is always well-defined for any $n>r>0$. As will be shown by Theorem  \ref{theorem1-epenalty} later, problem \eqref{eq:penprob} precisely serves as a global exact penalty of \eqref{eq:ori} provided that the feasible set of \eqref{eq:sh} is nonempty.
 \subsection{Our Contributions}

 The contribution of this work is threefold. Firstly, we introduce an equivalent model \eqref{eq:ori} that shares the same global and local optimal solutions as the original model \eqref{eq:sh}. By investigating the global exact penalty for model \eqref{eq:ori}, we establish a well-defined continuous optimization model \eqref{eq:penprob}. This continuous optimization model is equivalent to \eqref{eq:sh} in a global sense when $\rho$ surpasses a threshold $\overline{\rho}>0$ and the original model \eqref{eq:sh} is well-defined. In particular, it simultaneously takes into accout the binary constraint and the equality constraint in \eqref{eq:sh}. It is probably the main reason that our algorithm's the numerical performance surpass others.

 Secondly, we verify that the restricted Stiefel manifold $\bnr$ is actually an embedded submanifold of $\mathbb{R}^{n\times r}$. This means that problem \eqref{eq:penprob} is a composite optimization problem over the  submanifold $\bnr$. In order to leverage the benefits of smoothness in optimization, we replace $h_{\frac{1}{\sqrt{n}}}$ with its Moreau envelope env$_{\gamma h_{\frac{1}{\sqrt{n}}}}$ associated to parameter $\gamma>0$, resulting in a smooth manifold optimization problem
\begin{equation}\label{eq:sep1}
 \underset{X\in\bnr}{\min}\Theta_{\rho,\gamma}(X):=f(X) + \rho{\rm env}_{\gamma h_{\frac{1}{\sqrt{n}}}}(X).
\end{equation}
 Under a mild condition, model \eqref{eq:sep1} is also shown to be a global exact penalty of model \eqref{eq:ori} when $\rho$ surpasses a threshold $\overline{\rho}>0$ and  is well-defined.

 Thirdly, we develop a Riemannian gradient descent method with a line-search strategy to solve the manifold optimization problem \eqref{eq:sep1}. Leveraging the advantageous manifold structure and the exact penalty property of model  \eqref{eq:sep1}, the proposed algorithm yields the solutions that are closer to the feasible set of model \eqref{eq:sh} than those yielded by the existing solvers (see Section \ref{sec:feap}). The numerical experiments, conducted on both synthetic datasets and real-world datasets such as MNIST and CIFAR-10, demonstrate that our algorithm  significantly outperforms others in terms of feasibility violation. Moreover, our algorithm exhibits comparable or even superior performance to existing methods in terms of evaluation metrics related to the Hamming distance.

\subsection{Related Work}\label{1partB}

The two predominant approaches for tackling the discreteness of the feasible set $\mathcal{M}(n, r)$ in solving problem \eqref{eq:sh} are the relaxation and alternating maximization methods. There are numerous studies utilizing at least one of them. The relaxation method involves removing the binary constraint in \eqref{eq:sh}, transforming the problem into a continuous optimization one subjecting to equality constraints $X^{\top}X=nI_{r}$ and $X^{\top}v=0$. It then projects the resulting solution back to the binary space\cite{liu2011hashing, weiss2008spectral, weiss2012multidimensional, liu2014collaborative}. Despite its capacity to simplify the optimization process, it is evident that the projected solution often deviates from satisfying the equality constraints, leading to suboptimal outcomes.

Alongside relaxation, the alternating maximization method is frequently employed. This method introduces auxiliary variables to segregate the binary constraint from others, followed by solving a series of subproblems alternately \cite{gong2012iterative, yang2015robust, do2016learning, hu2017hashing, liu2020reinforced, chen2020strongly, qin2022joint}. However, solutions generated by this class of methods may fall short in adequately satisfying the feasibility conditions of \eqref{eq:sh}, primarily due to the inherent lack of convergence guarantees when addressing the binary orthogonal problem.

Beyond the aforementioned two methods, several other approaches have been explored. Shen et al. \cite{shen2016fast} introduced the quadratic penalty function for equality constraints $X^{\top}X\!=nI_{r}$ and $X^{\top}v=0$, employing a projected gradient method to solve the penalized binary optimization problem. Hoang et al. \cite{hoang2020unsupervised} also utilized the quadratic penalty function for the equality constraints but employed MPEC-EPM, a binary optimization technique, to solve the penalized problem. Xiong et al. \cite{xiong2021generalized} introduced two auxiliary equality constraints and employed ADMM to solve the quadratic penalty problem induced by these auxiliary equality constraints. However, due to the absence of leveraging the manifold structure, the solutions obtained from these methods may not conform well to the orthogonal constraint $X^{\top}X=nI_{r}$.

Motivated by the limitations of existing approaches, in this work, we propose a relaxation method by a novel continuous optimization reformulation for \eqref{eq:sh}, which is always well-defined. Specifically, we introduce an equivalent model \eqref{eq:ori}, and investigate its global exact penalty induced by the $\ell_1$ distance from the box set $\Lambda$ and its corresponding Moreau envelope.
 Considering that the penalty problem \eqref{eq:sep1} induced by the Moreau envelope is a smooth optimization over an embedded submaifold with a favorable structure, we develop a Riemannian gradient descent algorithm with a  line-search strategy to solve it, rather than addressing the nonsmooth penalty problem \eqref{eq:penprob}. The proposed method effectively leverages the manifold structure of $\bnr$ and the exact penalty property of ${\rm env}_{\gamma h_{\frac{1}{\sqrt{n}}}}$, ensuring a binay code with a low level of feasibility violation. Numerical comparisons with relaxation, alternating maximization, discrete proximal linearized minimization (DPLM) \cite{shen2016fast}, Deep Cross-modality Spectral Hashing(DCSH), and alternating binary matrix optimization (ABMO) \cite{xiong2021generalized} demonstrate the superior performance of our method in terms of feasibility violations. Furthermore, it achieves comparable or even superior performance to existing methods evaluated by the indices related to the Hamming distance.

\subsection{Organization}

The rest of this paper is organized as follows.
Section \ref{sec:pre} provides the preliminaries.
In Section \ref{sec:algo}, we establish that \eqref{eq:penprob} is a global exact penalty of \eqref{eq:ori} by deducing a global error bound for the feasible set of \eqref{eq:ori}. Section \ref{sec:algpg} develops a Riemannian gradient descent method, and gives its convergence results. In Section \ref{sec:exp}, we assess the numerical performance of our method through experiments conducted on both synthetic datasets and real-world datasets.
Finally, in Section \ref{sec:conc}, we provide concluding remarks.

\section{Preliminaries}\label{sec:pre}

As we are going to utilize the geometric properties of the Riemannian manifold, in this section, we provide the relevant concepts and theories that are necessary for comprehending the our research. In Riemannain geometry, a manifold is defined as a set of points that endowed with a locally Euclidean structure near each point. Let ${\cal M}$ represent a general Riemannian manifold. Given a point $X$ on ${\cal M}$, a tangent vector at $X$ is defined as the vector that is tangent to any smooth curves on ${\cal M}$ through $X$. For a given $X\in {\cal M}$, $T_{X}{\cal M}$ and $N_{X}{\cal M}$ denote the tangent and normal space to ${\cal M}$ at $X$, respectively. If the tangent spaces of a manifold are equipped with a smoothly varying inner product, the manifold is called Riemannian manifold\cite{petersen2006riemannian}. A Riemannian submanifold of $\mathbb{R}^{n\times r}$ is a submanifold equipped with a Riemannian metric that is induced from $\mathbb{R}^{n\times r}$.
\begin{figure}[H] 
\centering
\includegraphics[width=0.48\textwidth]{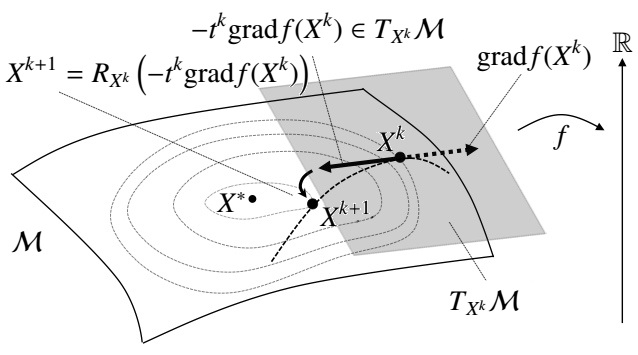}
\caption{\small Riemannian gradient descent}
    \begin{minipage}{0.45\textwidth}
        \small 
        \emph{Note:} This figure illustrates the iteration of search points in the Riemannian gradient descent method. It originates from \cite{HiroyukiSato2018}.
    \end{minipage}
\label{fig:cifar}
\end{figure}

The algorithm we employ belongs to the category of Riemannian gradient descent algorithm, which needs to compute the Riemannian gradient of the objective. For a compact submanifold $M$ (can be $\snr$ or $\bnr$, and the compactness of $\bnr$ is proved in Lemma \ref{lemma3:bnr}) of $\mathbb{R}^{n\times r}$, if a function  $h\!:\mathbb{R}^{n\times r}\to\overline{\mathbb{R}}\!:=\mathbb{R}\cup\{+\infty\}$ is differentiable at $X\in M$, then ${\rm grad}\,h(X)$ denotes the Riemannian gradient of $h$ at $X$. It is the unique tangent vector at $X$ such that for any $H\in{\rm T}_{X}M$, $\langle\nabla h(X),H\rangle=\langle{\rm grad}\,h(X),H\rangle$. For the Stiefel manifold, the tangent space $T_{X}\snr$ and the normal space $N_{X}\snr$ are given by
 \[
  T_{X}\snr:=\big\{H\in\R^{n\times r}\,|\,X\tp H +H\tp X = 0\big\},
 \]
 and
 \[
 N_{X}\snr:=\big\{XG\,|\,G\in\mathbb{S}^r\big\},
 \]
where $\mathbb{S}^r$ denotes the set of all $r\times r$ real symmetric matrices. The tangent spaces of ${\rm St}(n,r)$ are endowed with the trace inner product, i.e., $\langle \xi,\eta\rangle={\rm tr}( \xi\tp \eta)$ for any $\xi,\eta\in T_{X}\snr$, and its induced Frobenius norm $\|\cdot\|_F$.
For a closed set $\Omega\subseteq\Rnr$,  the notation $\delta_\Omega$ denotes the indicator function of $\Omega$,
 i.e., $\delta_\Omega(x)=0$ if $x\in \Omega$ and $+\infty$ otherwise, ${\rm proj}_{\Omega}$ denotes the projection operator onto $\Omega$ on the Frobenius norm, and ${\rm dist}(x,\Omega)\!:=\inf_{z\in\Omega}\|z\!-\!x\|_F=\|{\rm proj}_{\Omega}(x)-x\|_F$. The projection operator onto the tangent space of $\snr$ has the following expression
 \begin{equation}\label{eq:projst}
   {\rm proj}_{T_{X}\snr}(Z):=Z- X{\rm sym}(X\tp Z),
 \end{equation}
 where sym$(A) := (A + A\tp)/2$. Then, the Riemannian gradient of a differentiable function $h$ can be calculated by:
\[
	{\rm grad}h(X) = {\rm proj}_{T_{X}\snr}\lf(\nabla h(X)\rg).
\]

 Moving a point in a tangent vector direction while restricting on the manifold is realized by a retraction mapping\cite{Absil08}, whose formal definition is stated as follows.

\begin{definition}[retraction]\label{def:retr}
 A retraction on a manifold $\mathcal{M}$ is a $C^{\infty}$(smooth)-mapping $R$ from
   the tangent bundle\footnote{The domain of $R_X$ does not need to be the entire tangent bundle. However, the retraction we selected in this paper is the case.} $T \mathcal{M}:=\bigcup_{X\in\mathcal{M}}T_{X}\mathcal{M}$
   onto $\mathcal{M}$ and for any $X\in\mathcal{M}$, the restriction of $R$ to $T_{X}\mathcal{M}$,
   denoted by $R_X$, satisfies the following two conditions:
   \begin{itemize}
     \item [(i)] $R_{X}(0_{X})=X$, where $0_{X}$ denotes
                 the zero element of $T_{X}\mathcal{M}$;

     \item[(ii)] the differential of $R_X$ at $0_{X}$, ${\rm D}R_{X}(0_{X})$, is the identity mapping on $T_{X}\mathcal{M}$.
 \end{itemize}
\end{definition}
\begin{remark}\label{rmk:retr}
 For the submanifold $M$ of $\mathbb{R}^{n\times r}$,
  Item (ii) in Definition \ref{def:retr} is equivalent to saying that
  \[
   \lim_{{{\rm T}_{\!x}M}\ni\xi\to 0}\frac{\|R_{x}(\xi)-(x+\xi)\|_F}{\|\xi\|_F}=0.
  \]
 Common retractions for the Stiefel manifold include the exponential mapping, the Cayley transformation, the polar decomposition and the QR decomposition. The retractions defined on a compact submanifold often explicit a Lipschitz-type property; see the following lemma. 
 \end{remark}

\begin{lemma}\label{lem:pm}
  (see \cite[Eq. (B.3)$\,\&$(B.4)]{Boumal18})
  Let $M$ be a compact submanifold of a Euclidean space $\Rnr$,
  and let $R\!:{\rm T}M\to M$ be a retraction on $M$.
  Then, there exist constants $c_1>0$ and $c_2>0$ such that for any $X\in M$
  and $\xi\in {\rm T}_{\!X}M$,
  \[
   \|R_{X}(\xi)-X\|_F\le c_1\|\xi\|_F,
  \]
  and
  \[
  \|R_{X}(\xi)-(X+\xi)\|_F\le c_2\|\xi\|_F^2.
  \]
\end{lemma}

\subsection{Manifold Structure of $\bnr$}\label{sec2.1}

The set $\mathcal{B}_{v}(n,r)$ involved in \eqref{eq:ori} is an intersection of two compact embedded submanifolds of $\mathbb{R}^{n\times r}$: $\mathbb{L}_{v}$ and  ${\rm St}(n,r)$. Generally, the intersection of two embedded submanifolds of $\mathbb{R}^{n\times r}$ does  not necessarily result in an embedded submanifold. However, in our case, we can justify that $\mathcal{B}_{v}(n,r)$   is still a compact embedded submanifold of $\mathbb{R}^{n\times r}$.
 \begin{lemma}\label{lemma3:bnr}
 The set $\bnr$ is an embedded submanifold of $\Rnr$ and $\bnr$ is compact.
 \end{lemma}

By combining Lemma \ref{lemma3:bnr} with \cite[Exercise 6.7]{RW98}, we have the following characterization for the tangent and normal spaces of the manifold $\bnr$ at any $X\in\bnr$.
 \begin{lemma}\label{lemma1-tangent}
 Consider any $X\in\bnr$. The tangent and normal spaces of $\bnr$ at $X$ take the following form
 \begin{align*}
  T_X\bnr&= \big\{Z\in\Rnr\,|\, X\tp Z + Z\tp X = 0,\,Z\tp v = 0\big\},\\
  N_X\bnr&=\big\{XA + v\xi\tp\,|\,A\in\mathbb{S}^r,\,\xi\in\R^{r}\big\}.
 \end{align*}
\end{lemma}
 To characterize the projection mapping onto the tangent space $T_X\bnr$, we need the following lemma.
 \begin{lemma}\label{lemma2-tangent}
  Let $\mathbb{L}_1\!:=\{X\in\mathbb{R}^{n\times r}\,|\,\mathcal{A}_1(X)=0\}$ and $\mathbb{L}_2\!:=\{X\in\mathbb{R}^{n\times r}\,|\,\mathcal{A}_2(X)=0\}$ be two subspaces, where $\mathcal{A}_1\!:\mathbb{R}^{n\times r}\to\mathbb{R}^{n\times r}$ and $\mathcal{A}_2\!:\mathbb{R}^{n\times r}\to\mathbb{R}^r$ are the given linear mappings. If $\mathcal{A}_1\circ\mathcal{A}_2^*=0$ and $\mathcal{A}_2\circ\mathcal{A}_1^*=0$, then ${\rm proj}_{\mathbb{L}_1\cap\mathbb{L}_2}={\rm proj}_{\mathbb{L}_1}\circ{\rm proj}_{\mathbb{L}_2}={\rm proj}_{\mathbb{L}_2}\circ{\rm proj}_{\mathbb{L}_1}$ where, for any $Z\in\mathbb{R}^{n\times r}$,
  \begin{align*}
   {\rm proj}_{\mathbb{L}_1}(Z)=\big(\mathcal{I}-\mathcal{A}_1^*(\mathcal{A}_1\mathcal{A}_1^*)^{\dagger}\mathcal{A}_1\big)(Z),\\
   {\rm proj}_{\mathbb{L}_2}(Z)=\big(\mathcal{I}-\mathcal{A}_2^*(\mathcal{A}_2\mathcal{A}_2^*)^{\dagger}\mathcal{A}_2\big)(Z),
  \end{align*}
where the notation $\mathcal{Q}^{\dagger}$ denotes the pseudo-inverse of a positive semidefinite linear operator $\mathcal{Q}\!:\mathbb{R}^{n\times r}\to\mathbb{R}^{n\times r}$.
  \end{lemma}

Note that the tangent space of an embedded submanifold of $\mathbb{R}^{n\times r}$ is still an embedded submanifold. The tangent space $T_X\mathbb{R}^{n\times r}$ is the intersection of two embedded submanifolds: $T_X {\rm St}(n,r)$ and $\mathbb{L}_v$. These are associated with the linear mappings $\mathcal{A}_1:\mathbb{R}^{n\times r} \to \mathbb{S}^r$ and $\mathcal{A}_2:\mathbb{R}^{n\times r} \to \mathbb{R}^r$, defined by $\mathcal{A}_1(Z) := X^{\top}Z + Z^{\top}X$ and $\mathcal{A}_2(Z) := Z^{\top}v$, respectively. By Lemma \ref{lemma2-tangent}, the following result holds.
 \begin{corollary}\label{coro:pjT}
  Consider any $X\in\bnr$. The projection mapping onto $T_X\bnr$ is given by
  \[
    {\rm proj}_{T_X\bnr}(Z)={\rm proj}_{\mathbb{L}_{v}}(Z) - X{\rm sym}(X\tp{\rm proj}_{\mathbb{L}_{v}}(Z)),
  \]
  where ${\rm proj}_{\mathbb{L}_{v}}(Z):=\big(I_n - \frac{vv\tp}{\|v\|_2^2}\big)Z$ is the projection mapping onto  $\mathbb{L}_{v}$.
 \end{corollary}

It can be shown that many frequently used retractions onto the Stiefel manifold satisfy Definition \ref{def:retr} in the context of $\bnr$, such as the QR decomposition, exponential map, polar decomposition, and the Cayley transformation. Consequently, we follow the practice in \cite{chen2020proximal} and use the QR decomposition as the retraction onto the $\mathbb{R}^{n \times r}$ manifold in the numerical experiment. Specifically,
\[
	R_X(V) = {\rm qf}(X + V).
\]

\subsection{The Moreau Envelope of Function $h_{\frac{1}{\sqrt{n}}}$}\label{sec2.2}

To take advantage of smoothness in optimization, we need both the proximal mapping and the Moreau envelope associated with the penalty function $h_{\frac{1}{\sqrt{n}}}$. For a proper lower semicontinuous  function $h\!:\mathbb{R}^{n\times r}\to\overline{\mathbb{R}}$ and a constant $\gamma>0$, ${\rm prox}_{\!\gamma h}$ and ${\rm env}_{\gamma h}$ denote the proximal mapping and Moreau envelope of $h$ associated to $\gamma$, respectively, defined by
 \[
 {\rm prox}_{\!\gamma h}(X)
 \!:=\mathop{\arg\min}_{Z\in\mathbb{R}^{n\times r}}
 \Big\{\frac{1}{2\gamma}\|Z-X\|_F^2+h(Z)\Big\},
 \]
 and
 \[
  {\rm env}_{\gamma h}(X)\!:=\!\min_{Z\in\mathbb{R}^{n\times r}}
 \Big\{\frac{1}{2\gamma}\|Z-X\|_F^2+h(Z)\Big\}.
 \]
 When $h$ is convex, ${\rm prox}_{\!\gamma h}$ is a Lipschitz continuous mapping of modulus $1$ from $\Rnr$ to $\Rnr$, and $e_{\gamma h}$ is a continuously differentiable convex function with
 \(
  \nabla{\rm env}_{\gamma h}(X)=\gamma^{-1}(X\!-\!{\rm prox}_{\!\gamma h}(X)).
 \)
 The following lemma provides the proximal operator and the Moreau envelope of $h_{c}$,
 the $\ell_1$-distance function from the matrix box set $[-cE,cE]$ for a given $c>0$.
 \begin{lemma}\label{lem:prox}
  For $h_c(X):=\lng E, \max(0, X-c)-\min(0,X+c)\rng$ for $X\in\mathbb{R}^{n\times r}$, where $c>0$ is a given constant, its proximal operator and  Moreau envelope are respectively
  \[
   {\rm prox}_{\gamma h_c}(X)= \lf\{\begin{array}{cl}
			X_{ij}  & {\rm if}\;  |X_{ij}| < c, \\
			c\sgn(X_{ij})  &  {\rm if}\; c \le |X_{ij}|\le c + \gamma, \\
			\;X_{ij}-\gamma\sgn(X_{ij}) &{\rm if}\; |X_{ij}|> c + \gamma,
		\end{array}\rg.
  \]
  and ${\rm env}_{\gamma h}(X)=\sum_{i=1}^n\sum_{j=1}^r\theta_{\gamma,c}(X_{ij})$ with
  \[
  \theta_{\gamma,c}(x):=
   \lf\{\begin{array}{cl}
			0  & {\rm if}\;  |x| < c, \\
   \frac{1}{2\gamma}\lf(c\sgn(x) - x\rg)^2  &  {\rm if}\; c \le |x|\le c + \gamma, \\
		\;|x|-\frac{\gamma}{2} -c &{\rm if}\; |x|> c + \gamma.
	\end{array}\rg.
  \]
 \end{lemma}
\subsection{The stationary point of problem \eqref{eq:penprob}}
We recall from \cite{RW98} the basic (limiting/or Morduhovich) subdifferential of a generated real-valued function $h\!:\mathbb{X}\to\overline{\mathbb{R}}$.
 \begin{definition}\label{Gsubdiff-def}
 (see \cite[Definition 8.3]{RW98}) Consider a function $h\!:\mathbb{X}\to\overline{\mathbb{R}}$ and a point $x\in{\rm dom}\,h$. The regular subdifferential of $h$ at $x$ is defined as
  \[
   \widehat{\partial}h(x):=\Big\{v\in\mathbb{X}\,\lf|\,
   \liminf_{x\ne x'\to x}\frac{h(x')-h(x)-\langle v,x'-x\rangle}{\|x'-x\|}\ge 0\rg.\Big\},
  \]
  and the basic (also known as the limiting or Morduhovich) subdifferential of $h$ at $x$ is
  \[
   \partial h(x):=\lf\{v\in\mathbb{X}\lf\lvert\begin{array}{l} \exists\,x^k\to x\ {\rm with}\ h(x^k)\to h(x)\\
   {\rm and}\ v^k\in\widehat{\partial}h(x^k)\ {\rm s.t.}\ v^k\to v
\end{array}\rg.\rg\}.
  \]
 \end{definition}

Inspired by this, we introduce the following notion of stationary point.
\begin{definition}\label{def:sp}
	A matrix $X\in \bnr$ is called a stationary point of \eqref{eq:penprob} whenever
	\[
		0 \in {\rm proj}_{T_X\bnr}\lf(\nabla f(X)+ \rho \partial h_{\frac{1}{\sqrt{n}}}(X)\rg).
	\]
\subsection{The equivalence of models \eqref{eq:sh} and \eqref{eq:ori}}\label{sec2.3}

The equivalence between models \eqref{eq:sh} and \eqref{eq:ori} is demonstrated by the following lemma. The proof is presented in  \ref{app:equiv}.
 \begin{lemma}\label{equiv-model}
  Let $X=\frac{1}{\sqrt{n}}B$. Problem \eqref{eq:ori} is equivalent to \eqref{eq:sh} in the sense that they have the same global and local optimal solution sets. Moreover, each feasible point of \eqref{eq:ori} is a local optimal solution.
 \end{lemma}

\section{The Global Exactness of Problems \eqref{eq:penprob} and  \eqref{eq:sep1}}\label{sec:algo}

We will
 first justify that if model \eqref{eq:sh} has a nonempty feasible set, problem \eqref{eq:penprob} is a global exact penalty of \eqref{eq:sh}. In other words, there exists a threshold $\overline{\rho}>0$ such that the penalty problem \eqref{eq:penprob} associated to every $\rho\ge\overline{\rho}$ shares the same global optimal solution set as problem \eqref{eq:ori} does. To this end, we assume that $\mathcal{M}(n,r)\cap\mathbb{L}_{v}\ne\emptyset$ in this section.

 The following lemma shows that the distance of each point in $\Lambda$ from the set ${\rm St}(n,r)\cap\Lambda$ is locally bounded by its distance from $\snr$.
 \begin{lemma}\label{lemma-boxStnr}
  Consider any $\overline{X}\in{\rm St}(n,r)\cap\Lambda$. Then there exists a constant $\delta>0$ such that for all $X\in\mathbb{B}(\overline{X},\delta)\cap\Lambda$,
 \begin{equation*}
  {\rm dist}\big(X,{\rm St}(n,r)\cap\Lambda\big)\le 2\sqrt{nr}{\rm dist}(X,\snr).
 \end{equation*}
 \end{lemma}

 By invoking Lemma \ref{lemma-boxStnr}, we can derive a local Lipschitzian error bound for the constraint system $X\in\bnr\cap\Lambda$.
 \begin{proposition}\label{prop-BnrOmega}
  Consider any $\overline{X}\in\bnr\cap\Lambda$. Then there exists a constant $\delta>0$ such that for all $X\in\mathbb{B}(\overline{X},\delta)$,
  \begin{align*}
   &{\rm dist}\big(X,\bnr\cap\Lambda\big) \\
   \le&(2\sqrt{nr}\!+\!1)\big[{\rm dist}(X,\bnr)+{\rm dist}(X,\Lambda)\big].
  \end{align*}
 \end{proposition}
 \begin{remark}
  Proposition \ref{prop-BnrOmega} demonstrates that the metric qualification condition (see \cite[section 3.1]{Ioffe08}) holds for the set $\bnr\cap\Lambda$, which is equivalent to the metric subregularity of the multifunction
  $\mathcal{F}(X)\!:=(X;X)-\Lambda\times\bnr$ for $X\in\mathbb{R}^{n\times r}$ at $\overline{X}$ for the origin with the subregularity modulus $2\sqrt{nr}\!+\!1$. Although this result can be achieved by using \cite[Corollary 3.1]{BaiYe19} and  $\mathcal{T}_{\Lambda}(\overline{X})\cap\mathcal{T}_{\bnr}(\overline{X})=\{0\}$, the exact  subregularity modulus is unavailable in this way. Here, $\mathcal{T}_{\Lambda}(\overline{X})\cap\mathcal{T}_{\bnr}(\overline{X})=\{0\}$ is implied  by using $\overline{X}_{ij}\ne 0$ for all $(i,j)\in[n]\times[r]$ and combining the expression of $\mathcal{T}_{\bnr}(\overline{X})$ and that of $\mathcal{T}_{\Lambda}(\overline{X})$ given by
  \begin{align*}
  \mathcal{T}_{\Lambda}(\overline{X})
  &=\big\{H\in\mathbb{R}^{n\times r}\,|\,H_{ij}\ge 0\ {\rm for}\ (i,j)\in I_{-};\\
  &\qquad\qquad\qquad\quad H_{ij}\le0\ {\rm for}\ (i,j)\in I_{+}\big\},
  \end{align*}
  where $I_{+}\!:=\big\{(i,j)\in[n]\times[r]\ |\  \overline{X}_{ij}=\frac{1}{\sqrt{n}}\big\}$ and $I_{-}\!:=\big\{(i,j)\in[n]\times[r]\ |\ \overline{X}_{ij}=-\frac{1}{\sqrt{n}}\big\}$.
 \end{remark}

 By combining Proposition \ref{prop-BnrOmega} with the compactness of $\bnr$, the following global error bound holds.
\begin{corollary}\label{coro:gerror}
 There exists a constant $\kappa\!>0$ such that
 \[
  {\rm dist}\big(Z,\bnr\cap\Lambda\big)\le\kappa{\rm dist}(Z,\Lambda)
  \ \ \forall Z\in\bnr.
 \]
\end{corollary}

 Now we are prepared  to establish the global and local exact penalty results for problem \eqref{eq:penprob}. The proof of the exactness is inspired by \cite{qian2021exact}. Recalling that $f$ is continuously differentiable and the manifold $\bnr$ is compact, we note that $f$ is Lipschitz continuous on $\bnr$. Inspired by
 \begin{theorem}\label{theorem1-epenalty}
  Let $L_{\!f}$ be the Lipschitz modulus of $f$ on $\bnr$, and $f^*$ be the optimal value of \eqref{eq:ori}. Then, for all $X\in\bnr$,
  $f(X)-f^*+\kappa L_{\!f}\,h_{\!\frac{1}{\sqrt{n}}}(X)\ge 0$, where $\kappa$ is same as the one in Corollary \ref{coro:gerror}. Hence,  there exists a threshold $\overline{\rho}=\kappa L_{\!f}$ such that problem \eqref{eq:penprob} associated to each $\rho\ge\overline{\rho}$ shares  the same global optimal solution set as problem \eqref{eq:ori} does.
 \end{theorem}
 \begin{corollary}\label{coro:epenalty}
  For each local minimizer $X^*$ of \eqref{eq:ori}, there exist $\kappa'>0$ and  $\delta>0$ such that for all $X\in\mathbb{B}(X^*,\delta)\cap\bnr$,
  \[
	f(X)-f(X^*)+\kappa'L_{\!f}\,h_{\!\frac{1}{\sqrt{n}}}(Y)\ge 0,
  \]

Hence, each local minimizer of \eqref{eq:ori} is locally optimal to \eqref{eq:penprob} with $\rho\ge\!\kappa'L_{\!f}$. Conversely, each local minimizer of \eqref{eq:penprob} with $\rho>0$ in $\Lambda$ is locally optimal to \eqref{eq:ori}.
 \end{corollary}

 For problem \eqref{eq:sep1}, under a mild condition on $f$, we can justify that it is also an exact penalty problem of \eqref{eq:ori}.
 \begin{theorem}\label{thm:env_ep}
  Fix any $\gamma >0$. Suppose that for every global (or local) optimal solution $X^*$ of problem \eqref{eq:ori}, there exist $L'>0$ and $\delta'>0$ such that for all $X\in \mathbb{B}(X^*,\delta')\cap\bnr$ and all $\overline{X}\in{\rm proj}_{\bnr\cap\Lambda}(X)$,
  \begin{equation}\label{eq:epthm0}
  	f(X)\!-\!f(\overline{X}) \ge -L'\|X\!-\!\overline{X}\|_F^2.
  \end{equation}
  Then, for every global (or local) optimal solution $X^*$ of \eqref{eq:ori}, there exists $\varepsilon>0$ such that for all $X\in\mathbb{B}(X^*,\varepsilon)\cap\bnr$,
  \[
  	f(X)-f(X^*) + 2L'\gamma(2\sqrt{nr}\!+\!1)^2{\rm env}_{\gamma h_{\frac{1}{\sqrt{n}}}}(X)\ge 0,
  \]
  and consequently there exists a threshold $ \overline{\rho}>0$ such that problem \eqref{eq:sep1} associated to each $\rho\ge\overline{\rho}$ has the same global optimal solution set as problem \eqref{eq:ori} does. 
\end{theorem}

\section{Manifold Gradient Method for Problem \eqref{eq:sh}}\label{sec:algpg}

Building on
the results in Section \ref{sec:algo}, we know that problem \eqref{eq:penprob} associated to each $\rho \ge \kappa L_f$, as well as \eqref{eq:sep1} associated to each $\rho\ge2L'\gamma(2\sqrt{nr}\!+\!1)^2$ under the growth condition \eqref{eq:epthm0}, is equivalent to the discrete problem \eqref{eq:ori} or \eqref{eq:sh}. However, it is noteworthy that solving the nonsmooth penalty problem \eqref{eq:penprob} is more arduous than solving its smooth counterpart \eqref{eq:sep1}. We apply the manifold gradient descent method to seek an approximate stationary point of \eqref{eq:penprob} or \eqref{eq:sep1}.
\end{definition}

In each step, this method seeks a tangent direction $V^k$ to the manifold $\bnr$ at $X^k$, defined by
\begin{align}\label{subprob}
	V^k:=&\underset{V\in T_{X^k}\bnr}{\arg\min}\lf\{\lng \nabla \Theta_{\rho,\gamma}(X^k),V\rng+\frac{1}{2t_k}\|V\|_F^2\rg\}\nonumber \\
	=&-t_k\,{\rm grad}\,\Theta_{\rho,\gamma}(X^k),
\end{align}
and then use the retraction mapping $R_{X^k}$ to pull back this vector to $\bnr$. The efficiency of the manifold  gradient descent method heavily depends on the choice of step-size $t_k$, which in turn depends on the Lipschitz modulus $L_{\rho,\gamma}$ of $\nabla\Theta_{\rho,\gamma}$. When the Lipschitz modulus of $\nabla\!f$, denoted by $L_{\nabla\!f}$, is unavailable, $L_{\rho,\gamma}$ may be unknown. To address this issue, we develop the \textbf{M}anifold \textbf{G}radient method for \textbf{B}inary \textbf{O}rthogonal problem (MGBO), whose iterates are described in Algorithm \ref{ManPG1}.

\begin{algorithm}[h]
\caption{(Nonmonotone \textbf{M}anifold \textbf{G}radient method for \textbf{B}inary \textbf{O}rthogonal problem (MGBO))}\label{ManPG1}
\begin{algorithmic}[1]
\STATE Initialization: Select $\rho>0,\gamma>0$. Choose $\epsilon>0$, 
$\eta\in(0,1)$, $0<\alpha<1,0<t_{\rm min}\le t_{\rm max}$ and $X^0\in\bnr$.

\STATE \textbf{for}~{$k=0,1,2,\ldots$}

\STATE\label{Astep3}\quad Choose $t_k\in[t_{\rm min},t_{\rm max}]$;

\STATE \quad Let $V^k=-t_k\,{\rm grad}\,\Theta_{\!\rho,\gamma}(X^k)$
        and $X^{k+1}=R_{X^k}(V^k)$;
\STATE \quad\textbf{if} $\|{\rm grad}\,\Theta_{\!\rho,\gamma}(X^{k+1})\|\le \epsilon$, return;\label{algo:tc} \textbf{end if}.
\STATE\quad Let $J_k=\{\max\{0,k\!-\!m\},\ldots,k\}$;
\STATE\quad Set $\ell=0$, $t_{k,\ell} = t_k$, and $V^{k,\ell}=V^k$.
\STATE\quad\textbf{while} $\Theta_{\rho,\gamma}(X^{k+1})>{\displaystyle\max_{j\in J_{k}}}\,\Theta_{\rho,\gamma}(X^j)-\frac{\alpha}{2t_{k,\ell}}\|V^{k,\ell}\|_F^2$ \label{Astep5}\textbf{do}
\STATE\label{Astep6} \quad\quad Let $t_{k,\ell+1}=\eta t_{k,\ell}$;
\STATE \quad\quad Set  $V^{k,\ell+1}= -t_{k,\ell+1}\,{\rm grad}\,\Theta_{\!\rho,\gamma}(X^k)$;
\STATE\label{Astep7} \quad\quad Set $X^{k+1}=R_{X^k}(V^{k,\ell+1})$;
\STATE \quad\textbf{end while}
\STATE \textbf{end for}
\end{algorithmic}
\end{algorithm}

 A good step-size initialization at each outer iteration can greatly reduce the line-search cost. Inspired by \cite{Wright09}, we initialize the step-size $t_k$ by the Barzilai-Borwein rule \cite{Barzilar88}:
 \begin{equation}\label{BB-stepsize}
  t_k\!=\!\max\Big\{\!\min\Big\{\!\frac{\|\Delta\!X^k\|_F^2}{|\langle\Delta\!X^k,\Delta Y^k\rangle|},
  \frac{|\langle\Delta\!X^k,\Delta Y^k\rangle|}{\|\Delta Y^k\|_F^2},t_{\rm max}\Big\},
  t_{\rm min}\Big\}
 \end{equation}
 where $\Delta X^k\!:=X^k\!-\!X^{k-1}$ and $\Delta Y^k\!:={\rm grad}\,\Theta_{\rho,\gamma}(X^k) \!-\!{\rm grad}\,\Theta_{\rho,\gamma}(X^{k-1})$. The step size is accepted once the nonmonotone line-search criterion (refer to \cite{Grippo86,Grippo02}) in Step \ref{Astep5} is met. Additionally, by following the same proof as for \cite[Lemma 8]{qian2021exact}, the nonmonotone line-search criterion specified in Step \ref{Astep5} is well defined.
 \begin{lemma}\label{lemma1-ManPG1}
  For each $k\in\mathbb{N}$, the nonmonotone line-search criterion in Step \ref{Astep5} is satisfied whenever $t_k\le\bar{t}:=\frac{1-\alpha}{2c_2c_{\rho,\gamma}
 		+c_1^2L_{\rho,\gamma}}$, where $c_{\rho,\gamma}\!:={\displaystyle\max_{Z\in{\rm St}(n,r)}}
 	\|\nabla\Theta_{\rho,\gamma}(Z)\|_F$ and $c_1,c_2$ are the constants from Lemma \ref{lem:pm}.
 \end{lemma}

 By the expressions of $\Theta_{\rho,\gamma}$ and ${\rm grad}\,\Theta_{\rho,\gamma}$, we have that
 \begin{align*}
  &\frac{1}{\rho}\|{\rm grad}\,\Theta_{\rho,\gamma}(X^k)-{\rm grad} f(X^k)\|_F\\
  &=\big\|\gamma^{-1}{\rm proj}_{T_{X^k}\bnr}\big(X^k\!-{\rm prox}_{\gamma h_{\frac{1}{\sqrt{n}}}}(X^k)\big)\big\|_F,
 \end{align*}
which implies that for an appropriately large $\rho>0$, the projection of $X^k\!-{\rm prox}_{\gamma h_{\frac{1}{\sqrt{n}}}}(X^k)$ (the violation vector of $X^k$ from the box set $\Lambda$) onto the tangent space $T_{X^k}\bnr$ will be small. This means that Algorithm \ref{ManPG1} with an appropriately large $\rho>0$ will return a point $X^{k}\in\bnr$ with small $\big\|\gamma^{-1}{\rm proj}_{T_{X^k}\bnr}\big(X^k\!-{\rm prox}_{\gamma h_{\frac{1}{\sqrt{n}}}}(X^k)\big)\big\|_F$.


 For Algorithm \ref{ManPG1} with $\varepsilon=0$, we have the following global convergence result. Since its proof resembles to that of Lemma 4.6 and Theorem 4.7 in \cite{qian2021exact}, we here omit it.
 \begin{theorem}\label{theorem1-ManPG1}
  Let $\{(X^k,V^k)\}_{k=0}^{\infty}$ be the sequence generated by Algorithm \ref{ManPG1}. For each $k$, write $\ell(k):=\mathop{\arg\max}_{j\in J_k}
 	\Theta_{\rho,\gamma}(X^j)$. Then, the following assertions hold.
  \begin{enumerate}
   \item The sequences $\{\Theta_{\rho,\gamma}(X^k)\}_{k\in\mathbb{N}}$
 		and~~ $\{V^k\}_{k\in\mathbb{N}}$ are convergent, and $\lim_{k\to\infty}$ $V^k=0$.
 		
   \item The sequence $\{X^k\}_{k\in\mathbb{N}}$ is bounded and every cluster point
 		is a stationary point of \eqref{eq:sep1}, i.e., $\|{\rm grad}\,\Theta_{\!\rho,\gamma}(X)\|=0$ for $X\in\bnr$.
 		
   \item $\sum_{k=0}^{\infty}\|X^{k+1}\!-\!X^k\|_F<\!\infty$ provided that $f$ is definable in an o-minimal structure $\mathscr{O}$ over $\mathbb{R}$, and that
 	\begin{equation}\label{mild-cond}
 	 \sum_{\mathcal{K}\ni k=0}^{\infty}
 			\sqrt{\Theta_{\rho,\gamma}(X^{\ell(k+1)})-\Theta_{\rho,\gamma}(X^{k+1})}<\infty
 	\end{equation}
 	whenever ${\displaystyle
 		\liminf_{\mathcal{K}\ni k\to\infty}
 		\frac{\Theta_{\rho,\gamma}(X^{\ell(k)})-\Theta_{\rho,\gamma}(X^{\ell(k+1)})}{\|X^{k+1}-X^k\|_F^2}=0}$, where
 	 $\mathcal{K}\!:=\!\big\{k\in\mathbb{N}\ \lf|\
 		\Theta_{\rho,\gamma}(X^{\ell(k+1)})-\Theta_{\rho,\gamma}(X^{k+1})
 		\ge\frac{\alpha\|X^{k+1}-X^k\|_F^2}{6t_{\max}}\rg.\big\}$.
 	\end{enumerate}
 \end{theorem}
The following iteration complexity result of MGBO is inspired by \cite{Boumal18,chen2020proximal}.
\begin{theorem}\label{thm:itr}
Algorithm  \ref{ManPG1} returns $X\in\bnr$ satisfying $\|{\rm grad}\,\Theta_{\!\rho,\gamma}(X)\|_F\le \epsilon$ in at most
\[
\lf\lceil (m+1)(\frac{2(\Theta_{\rho,\gamma}(X^0) - \Theta_{\rho,\gamma}^*)}{(\alpha\eta \bar{t}\epsilon^2)}+1)\rg\rceil
\]
iterations, where $\bar{t}$ is defined in Lemma \ref{lemma1-ManPG1}, and $\Theta_{\rho,\gamma}^*$  is the optimum of \eqref{eq:sep1}.
\end{theorem}

 \section{Numerical Experiment}\label{sec:exp}
We 
 test the performance of MGBO to semantic hashing on both the randomly-generated and the real world datasets, and compare its performance with that of several existing solvers. All experiments are performed in MATLAB on a workstation running on 64-bit Windows System with an Intel(R) Core(TM) i9-12900 CPU 2.40 GHz and 64 GB RAM.

 Unless otherwise stated, the following parameters of Algorithm \ref{ManPG1} are used for the subsequent numerical tests:
 \begin{align*}
  \rho=10,\,\gamma=0.2,\,\epsilon = 10^{-5}\sqrt{n},m=5,\qquad\\
  \eta = 0.85,\,\alpha=10^{-4},\,t_{\min} = 10^{-20},\,t_{\max} = 10^{20}.
 \end{align*}
 The step-size $t_k$ in each outer iteration is initialized by the BB rule in \eqref{BB-stepsize}.

\begin{remark}
The choice of the $\rho$ and $\gamma$ has certain impact on the performance of the algorithm. For the penalty parameter $\rho$, empirical observations indicate that for $\rho < 10$, increasing its value correlates with a reduction in infeasibility. Interestingly, this trend stabilizes for $\rho \geq 10$, with imfeasibility remaining constant irrespective of further increases in $\rho$. The choice of $\gamma$ governs the trade-off between the amount of smoothing applied to penalty function $h_{\frac{1}{\sqrt{n}}}$ and the convergence speed. The rest of the hyperparameters are set as default in \cite{qian2021exact} and do not impact the algorithm's stability.
\end{remark}
\subsection{Numerical Performance in Finding a Feasible Point}\label{sec:feap}

 We test the efficiency of MGBO in finding a feasible point of problem \eqref{eq:sh} with  $v=e$, where
 $A$ is a Laplacian matrix generated randomly.  We conduct a comparative analysis of MGBO against the Relaxiation method (Relax.), alternating maximization method (Alt.Max.), alternating binary matrix optimization (ABMO) method \cite{xiong2021generalized}, the discrete proximal linearized minimization (DPLM) method \cite{shen2016fast}, and DCSH model \cite{hoang2020unsupervised} empowered by the MPEC-EPM. For the ABMO, we choose the LBFGS with Wolfe line-search rule to solve its subproblem and set all ABMO parameters to default values. For DPLM and DCSH method, we also use their default parameters. Since most of the above methods lack convergence guarantees, the maximum iteration count probably impact computation speed. Therefore, we set the maximum iteration for each method to \texttt{maxIter=1000}. The solutions obtained by all six solvers are rounded using the sign function, yielding the corresponding binary matrices.

To conduct numerical tests in finding a feasible point, we generate synthetic samples using the following parameters of dimension pairs:
\begin{subequations}
	\begin{equation}\label{eq:nt1}
	\{(n,t)~|~n = 2^a, r = a\},
	\end{equation}
	\begin{equation}\label{eq:nt2}
	\{(n,t)~|~n = 2^a, r = 2a\},
	\end{equation}
	\begin{equation}\label{eq:nt3}
	\{(n,t)|~n = 2^a, r = 2a+1\},
	\end{equation}
\end{subequations}	
for $a\in\{2,~3,~4,~5,~6,~7\}$. For example, when $a=2$ in \eqref{eq:nt1}, we have the dimension pair $n=2^2=4$ and $r=2$. The feasible set is ${\cal M}(4,2)$. As the value of $a$ increases, finding a feasible solution becomes progressively more difficult. Furthermore, we intensify the difficulty of finding a feasible solution by increasing the ratio of $r/n$ from $a/2^a$ to $(2a+1)/2^a$. It is worth noting that, as indicated by the results in \cite[Chapter 2]{hedayat1999orthogonal}, the feasible set of problem \eqref{eq:sh} is expected to be non-empty in these cases.

\begin{table}[H]
    \centering
    \small
    \caption{Feasibility satisfaction test on synthetic data}\label{tab:fea}
\begin{tabular}{llrrrrrr}
    \hline
\multicolumn{1}{l}{}    &         & \multicolumn{5}{c}{$a$}     &     \\
\multicolumn{1}{l}{}    &         & 2  & 3  & 4   & 5   & 6   & 7   \\\toprule
\multicolumn{1}{l}{$r=a$} &         &    &    &     &     &     &     \\
\multirow{2}{*}{MGBO}   & Balance & 0  & 0  & 0   & 26  & 80  & 98  \\
                        & Orth    & 0  & 0  & 6   & 48  & 92  & 100  \\
\multirow{2}{*}{Relax.} & Balance & 85 & 100 & 100 & 100 & 100 & 100 \\
                        & Orth    & 52 & 97 & 100 & 100 & 100 & 100 \\
\multirow{2}{*}{Alt.Max.}& Balance & 85 & 98 & 100 & 100 & 100 & 100 \\
                        & Orth    & 52 & 92 & 100 & 100 & 100 & 100 \\
\multirow{2}{*}{ABMO}   & Balance & 51 & 96 & 99  & 100 & 99 & 100 \\
                        & Orth    & 37 & 91 & 100 & 100 & 100 & 100 \\
\multirow{2}{*}{DPLM}   & Balance & 51 & 95 & 100 & 100 & 100 & 100 \\
                        & Orth    & 15 & 75 &  99 & 100 & 100 & 100 \\
\multirow{2}{*}{DCSH}& Balance & 57 & 100 & 100 & 100 & 100 & 100 \\
                        & Orth    & 19 & 93 & 100 & 100 & 100 & 100 \\\toprule
\multicolumn{1}{l}{$r=2a$}  &         & &     &     &     &     &     \\
\multirow{2}{*}{MGBO}   & Balance & - & 0   & 31  & 91  & 99  & 100  \\
                        & Orth    & - & 0   & 85  & 100  & 100 & 100 \\
\multirow{2}{*}{Relax.} & Balance & - & 100 & 100 & 100 & 100 & 100 \\
                        & Orth    & - & 100 & 100 & 100 & 100 & 100 \\
\multirow{2}{*}{Alt.Max.}& Balance & - & 100 & 100 & 100 & 100 & 100 \\
                        & Orth    & - & 100 & 100 & 100 & 100 & 100 \\
\multirow{2}{*}{ABMO}   & Balance & - & 99  & 100 & 100 & 100 & 100 \\
                        & Orth    & - & 99  & 100 & 100 & 100 & 100 \\
\multirow{2}{*}{DPLM}   & Balance & - & 100 & 100 & 100 & 100 & 100 \\
                        & Orth    & - & 100 & 100 & 100 & 100 & 100 \\
\multirow{2}{*}{DCSH}& Balance & - & 100 & 100 & 100 & 100 & 100 \\
                        & Orth    & - & 100 & 100 & 100 & 100 & 100 \\\toprule
\multicolumn{1}{l}{$r=2a\!+\!1$} &         & &     &     &     &     &     \\
\multirow{2}{*}{MGBO}   & Balance & - & 0   & 46  & 95  & 100 & 100 \\
                        & Orth    & - & 0   & 91  & 100 & 100 & 100 \\
\multirow{2}{*}{Relax.} & Balance & - & 100 & 100 & 100 & 100 & 100 \\
                        & Orth    & - & 100 & 100 & 100 & 100 & 100 \\
\multirow{2}{*}{Alt.Max.}& Balance & - & 100 & 100 & 100 & 100 & 100 \\
                        & Orth    & - & 100 & 100 & 100 & 100 & 100 \\
\multirow{2}{*}{ABMO}   & Balance & - & 100 & 100 & 100 & 100 & 100 \\
                        & Orth    & - & 100 & 100 & 100 & 100 & 100 \\
\multirow{2}{*}{DPLM}   & Balance & - & 100 & 100 & 100 & 100 & 100 \\
                        & Orth    & - & 100 & 100 & 100 & 100 & 100 \\
\multirow{2}{*}{DCSH}& Balance & - & 100 & 100 & 100 & 100 & 100 \\
                        & Orth    & - & 100 & 100 & 100 & 100 & 100 \\\bottomrule
\end{tabular}
    \begin{flushleft}
    {\small Note: The experiment for $a=2$ in $r=2a$ and $r=2a+1$ are ignored due to the fact that ${\cal M}(4,4)={\cal M}(4,5)=\emptyset$. }
    \end{flushleft}
\end{table}

 Table \ref{tab:fea} presents the counts of cases in which the binary matrices returned by three solvers violate either the orthogonal constraint $X\tp X=nI_r$ or the balanced constraint $X\tp e=0$, denoted by \texttt{Orth} or \texttt{Balance}, respectively, across 100 test instances. The results show that MGBO successfully yield solutions that satisfy the two constraints when $t=2$ and $t=3$. As the dimensions increase, the number of solutions violating the constraints also increases, but MGBO still demonstrates the highest likelihood of generating feasible solutions.


\subsection{Results for the Unsupervised Hashing}\label{sec:unha}

\begin{table*}[h]
\centering
\scriptsize
\caption{Numerical result for synthetic data: $r=16~\&~32$}\label{tab:exp_r}
\resizebox{\textwidth}{!}{\begin{tabular}{l|lllll|rrrrr|rrrrr}
\toprule
      & \multicolumn{5}{c}{CPU Time}                         & \multicolumn{5}{c}{Violation of   Constraint: Balance} & \multicolumn{5}{c}{Violation of   Constraint: Orthogonal} \\\hline
$r=16$ & & & & & & & & & & \\
$n=$         & $500$      & $1,000$    & $2,000$    & $5,000$    & $10,000$   & $500$       & $1,000$     & $2,000$    & $5,000$    & $10,000$   & $500$       & $1,000$     & $2,000$     & $5,000$     & $10,000$    \\\hline
\underline{MGBO} & $2.46^{-2}$ & $1.78^{-2}$ & $2.83^{-2}$ & $7.62^{-2}$ & $3.01^{-1}$ & \textbf{$1.60^{1}$} & \textbf{$4.09^{1}$} & \textbf{$3.52^{1}$} & \textbf{$4.63^{1}$} & \textbf{$1.34^{2}$} & \textbf{$9.78^{1}$} & \textbf{$1.59^{2}$} & \textbf{$1.96^{2}$} & \textbf{$3.79^{2}$} & \textbf{$7.84^{2}$} \\
\underline{Relax.} & \textbf{$1.78^{-2}$} & \textbf{$3.99^{-3}$} & \textbf{$6.53^{-3}$} & \textbf{$1.15^{-2}$} & \textbf{$2.65^{-2}$} & $3.67^{2}$ & $6.72^{2}$ & $1.62^{3}$ & $3.31^{3}$ & $5.78^{3}$ & $3.81^{2}$ & $5.59^{2}$ & $1.35^{3}$ & $2.19^{3}$ & $3.08^{3}$ \\
\underline{Alt.Max.} & $2.86^{-2}$ & $3.99^{-2}$ & $6.47^{-2}$ & $1.19^{0}$ & $3.35^{0}$ & $4.40^{1}$ & $8.47^{1}$ & $9.16^{1}$ & $1.37^{2}$ & $2.44^{2}$ & $2.87^{3}$ & $4.28^{3}$ & $2.69^{4}$ & $2.34^{4}$ & $3.17^{4}$\\
\underline{ABMO} & $6.00^{-1}$ & $2.23^{0}$ & $7.66^{0}$ & $3.87^{1}$ & $4.29^{2}$ & $3.83^{1}$ & $7.44^{1}$ & $1.49^{2}$ & $2.02^{2}$ & $2.05^{2}$ & $2.70^{2}$ & $3.74^{2}$ & $5.46^{2}$ & $9.14^{2}$ & $1.17^{3}$ \\
\underline{DPLM} & $1.21^{-1}$ & $1.63^{-1}$ & $6.41^{-1}$ & $1.23^{0}$ & $5.14^{0}$ & $6.35^{1}$ & $6.88^{1}$ & $4.80^{1}$ & $1.60^{2}$ & $1.04^{2}$ & $7.25^{3}$ & $1.32^{4}$ & $3.10^{4}$ & $7.75^{4}$ & $1.55^{5}$ \\
\underline{DCSH} & $3.24^{-2}$ & $3.77^{-2}$ & $1.63^{-1}$ & $2.40^{-1}$ & $1.03^{0}$ & $4.24^{2}$ & $4.00^{1}$ & $1.20^{2}$ & $1.04^{2}$ & $3.20^{1}$ & $7.75^{3}$ & $1.47^{4}$ & $3.10^{4}$ & $7.75^{4}$ & $1.55^{5}$ \\
\bottomrule
$r=32$ & & & & & & & & & & \\
$n=$         & $500$      & $1,000$    & $2,000$    & $5,000$    & $10,000$   & $500$       & $1,000$     & $2,000$    & $5,000$    & $10,000$   & $500$       & $1,000$     & $2,000$     & $5,000$     & $10,000$    \\\hline
\underline{MGBO} & $3.53^{-2}$ & $3.61^{-2}$ & $6.73^{-2}$ & $2.86^{-1}$ & $5.75^{-1}$ & \textbf{$2.02^{1}$} & \textbf{$2.79^{1}$} & \textbf{$6.32^{1}$} & \textbf{$9.12^{1}$} & \textbf{$5.18^{1}$} & \textbf{$1.78^{2}$} & \textbf{$3.14^{2}$} & \textbf{$4.72^{2}$} & \textbf{$7.58^{2}$} & \textbf{$1.59^{3}$} \\
\underline{Relax.}   & \textbf{$1.86^{-2}$} & \textbf{$5.71^{-3}$} & \textbf{$9.37^{-3}$} & \textbf{$2.57^{-2}$} & \textbf{$4.47^{-2}$} & $6.64^{2}$ & $1.17^{3}$ & $1.91^{3}$ & $4.69^{3}$ & $1.12^{4}$ & $1.02^{3}$ & $1.52^{3}$ & $2.05^{3}$ & $4.31^{3}$ & $1.22^{4}$ \\
\underline{Alt.Max.} & $3.45^{-2}$ & $6.58^{-2}$ & $1.58^{-1}$ & $2.23^{0}$ & $5.79^{0}$ & $7.93^{1}$ & $1.05^{2}$ & $1.47^{2}$ & $3.04^{2}$ & $2.75^{2}$ & $5.74^{3}$ & $8.97^{3}$ & $5.10^{4}$ & $7.23^{4}$ & $1.69^{5}$ \\
\underline{ABMO} & $2.51^{0}$ & $7.86^{0}$ & $2.56^{1}$ & $4.24^{2}$ & $1.47^{3}$ & $8.08^{1}$ & $1.06^{2}$ & $1.64^{2}$ & $2.29^{2}$ & $3.34^{2}$ & $5.16^{2}$ & $7.69^{2}$ & $1.08^{3}$ & $1.69^{3}$ & $2.40^{3}$ \\
\underline{DPLM} & $1.83^{-1}$ & $2.49^{-1}$ & $8.38^{-1}$ & $3.73^{0}$ & $7.16^{0}$ & $6.79^{1}$ & $1.52^{2}$ & $6.79^{1}$ & $1.02^{2}$ & $1.39^{2}$ & $1.57^{4}$ & $2.94^{4}$ & $6.30^{4}$ & $1.57^{5}$ & $3.14^{5}$ \\
\underline{DCSH} & $6.10^{-2}$ & $6.72^{-2}$ & $1.88^{-1}$ & $9.24^{-1}$ & $1.55^{0}$ & $2.26^{1}$ & $3.58^{1}$ & $6.00^{2}$ & $3.51^{2}$ & $3.17^{2}$ & $1.57^{4}$ & $3.07^{4}$ & $6.30^{4}$ & $1.57^{5}$ & $3.15^{5}$ \\
\bottomrule
\end{tabular}}
    \begin{flushleft}
    {\small Note: The superscript $k$ denotes a scale of $10^k$. }
    \end{flushleft}
\end{table*}

\begin{table*}[h]
\centering
\scriptsize
\caption{Numerical result for synthetic data: $n=1,000~\&~5,000$}\label{tab:exp_n}
\resizebox{\textwidth}{!}{\begin{tabular}{l|llllll|rrrrrr|rrrrrr}
\toprule
      & \multicolumn{6}{c}{CPU Time}                         & \multicolumn{6}{c}{Violation of   Constraint: Balance} & \multicolumn{6}{c}{Violation of   Constraint: Orthogonal} \\\hline
$r=$         & $16$ & $32$ & $48$ & $64$ & $80$ & $96$ & $16$ & $32$ & $48$ & $64$ & $80$ & $96$ & $16$ & $32$ & $48$ & $64$ & $80$ & $96$ \\
$n=1,000$ & & & & & & & & & & & & & & & & & & \\\hline
\underline{MGBO} & $2.39^{-2} $ & $3.11^{-2} $ & $4.46^{-2} $ & $4.59^{-2} $ & $6.61^{-2} $ & $7.68^{-2} $ & \textbf{$9.44^{1} $} & \textbf{$2.37^{2} $} & \textbf{$3.46^{2} $} & \textbf{$4.42^{2} $} & \textbf{$5.75^{2} $} & \textbf{$6.97^{2} $} & \textbf{$1.89^{1} $} & \textbf{$2.86^{1} $} & \textbf{$3.41^{1} $} & \textbf{$4.29^{1} $} & \textbf{$4.53^{1} $} & \textbf{$5.32^{1} $} \\
\underline{Relax.}  & \textbf{$5.88^{-3} $} & \textbf{$2.07^{-3} $} & \textbf{$2.49^{-3} $} & \textbf{$2.93^{-3} $} & \textbf{$2.67^{-2} $} & \textbf{$4.13^{-3} $} & $4.18^{3} $ & $8.02^{3} $ & $9.24^{3} $ & $1.01^{4} $ & $1.09^{4} $ & $1.17^{4} $ & $2.12^{3} $ & $2.66^{3} $ & $2.72^{3} $ & $2.73^{3} $ & $2.73^{3} $ & $2.74^{3} $ \\
\underline{Alt.Max.} & $5.77^{-3} $ & $9.83^{-3} $ & $1.20^{-2} $ & $1.56^{-2} $ & $2.10^{-2} $ & $2.44^{-2} $ & $5.58^{3} $ & $1.18^{4} $ & $1.82^{4} $ & $2.45^{4} $ & $3.08^{4} $ & $3.72^{4} $ & $1.15^{2} $ & $1.28^{2} $ & $1.84^{2} $ & $2.25^{2} $ & $2.24^{2} $ & $2.61^{2} $ \\
\underline{ABMO} & $1.63^{0} $ & $4.04^{0} $ & $7.36^{0} $ & $1.17^{1} $ & $1.79^{1} $ & $2.56^{1} $ & $3.81^{2} $ & $7.59^{2} $ & $1.16^{3} $ & $1.56^{3} $ & $1.92^{3} $ & $2.28^{3} $ & $1.01^{2} $ & $9.82^{1} $ & $1.40^{2} $ & $1.51^{2} $ & $1.56^{2} $ & $1.73^{2} $ \\
\underline{DPLM} & $2.10^{-3} $ & $1.98^{-3} $ & $2.84^{-3} $ & $3.36^{-3} $ & $4.18^{-3} $ & $5.08^{-3} $ & $1.55^{4} $ & $3.15^{4} $ & $4.75^{4} $ & $6.35^{4} $ & $7.95^{4} $ & $9.50^{4} $ & $1.44^{2} $ & $1.24^{2} $ & $3.33^{2} $ & $8.00^{1} $ & $5.37^{1} $ & $2.06^{2} $ \\
\underline{DCSH} & $2.56^{-3} $ & $2.22^{-3} $ & $3.19^{-3} $ & $3.51^{-3} $ & $4.99^{-3} $ & $5.19^{-3} $ & $1.55^{4} $ & $3.15^{4} $ & $4.70^{4} $ & $6.35^{4} $ & $7.95^{4} $ & $9.55^{4} $ & $3.20^{1} $ & $2.16^{3} $ & $1.80^{3} $ & $2.72^{2} $ & $3.40^{2} $ & $3.53^{2} $ \\
\bottomrule
$r=$         & $16$ & $32$ & $48$ & $64$ & $80$ & $96$ & $16$ & $32$ & $48$ & $64$ & $80$ & $96$ & $16$ & $32$ & $48$ & $64$ & $80$ & $96$ \\
$n=5,000$ & & & & & & & & & & & & & & & & & & \\\hline
\underline{MGBO} & $6.14^{-2} $ & $2.60^{-1} $ & $4.08^{-1} $ & $4.88^{-1} $ & $6.50^{-1} $ & $7.60^{-1} $ & \textbf{$4.71^{1} $} & \textbf{$7.41^{1} $} & \textbf{$7.89^{1} $} & \textbf{$8.20^{1} $} & \textbf{$1.09^{2} $} & \textbf{$9.38^{1} $} & \textbf{$2.29^{2} $} & \textbf{$4.49^{2} $} & \textbf{$6.61^{2} $} & \textbf{$9.28^{2} $} & \textbf{$1.24^{3} $} & \textbf{$1.51^{3} $} \\
\underline{Relax.} & \textbf{$6.12^{-3} $} & \textbf{$6.11^{-3} $} & \textbf{$8.94^{-3} $} & \textbf{$1.09^{-2} $} & \textbf{$1.22^{-2} $} & \textbf{$1.42^{-2} $} & $1.11^{4} $ & $1.42^{4} $ & $1.47^{4} $ & $1.48^{4} $ & $1.48^{4} $ & $1.48^{4} $ & $2.25^{4} $ & $4.43^{4} $ & $5.16^{4} $ & $5.60^{4} $ & $5.99^{4} $ & $6.34^{4} $ \\
\underline{Alt.Max.} & $1.12^{-2} $ & $6.41^{-2} $ & $8.36^{-2} $ & $1.17^{-1} $ & $1.47^{-1} $ & $1.55^{-1} $ & $1.65^{2} $ & $2.32^{2} $ & $3.07^{2} $ & $3.42^{2} $ & $4.28^{2} $ & $4.85^{2} $ & $1.90^{4} $ & $4.33^{4} $ & $6.82^{4} $ & $9.22^{4} $ & $1.16^{5} $ & $1.40^{5} $ \\
\underline{ABMO} & $2.10^{1} $ & $3.58^{2} $ & $7.32^{2} $ & $1.18^{3} $ & $1.90^{3} $ & $2.51^{3} $ & $1.43^{2} $ & $3.21^{2} $ & $2.84^{2} $ & $3.08^{2} $ & $3.70^{2} $ & $4.16^{2} $ & $8.19^{2} $ & $1.75^{3} $ & $2.64^{3} $ & $3.50^{3} $ & $4.36^{3} $ & $5.21^{3} $ \\
\underline{DPLM} & $3.97^{-3} $ & $2.91^{-2} $ & $3.71^{-2} $ & $5.29^{-2} $ & $6.57^{-2} $ & $7.09^{-2} $ & $2.80^{2} $ & $5.20^{2} $ & $6.37^{2} $ & $1.44^{2} $ & $2.10^{2} $ & $5.49^{2} $ & $7.75^{4} $ & $1.57^{5} $ & $2.37^{5} $ & $3.17^{5} $ & $3.89^{5} $ & $4.77^{5} $ \\
\underline{DCSH} & $5.17^{-3} $ & $3.87^{-2} $ & $5.04^{-2} $ & $6.75^{-2} $ & $8.29^{-2} $ & $8.82^{-2} $ & $8.00^{1} $ & $3.17^{2} $ & $6.24^{2} $ & $4.00^{2} $ & $2.73^{2} $ & $9.80^{1} $ & $7.75^{4} $ & $1.57^{5} $ & $2.37^{5} $ & $3.17^{5} $ & $2.97^{5} $ & $4.77^{5} $ \\
\bottomrule
\end{tabular}}
\end{table*}

 We next conduct tests for the unsupervised spectral hashing (SH) \cite{weiss2008spectral}, the model is described as:
\begin{align}\label{eq:unsup_dh}
 &\underset{B\in\{-1, 1\}^{n\times r}}{\min} \widetilde{f}(B) := {\rm tr}(B\tp A B) , \nonumber\\
 &\qquad\ {\rm s.t.~}\ B\tp B = nI_r,\quad B\tp e = 0,
\end{align}
where $A$ is a Laplacian matrix. Model \eqref{eq:unsup_dh} can be formulated as the equivalent model \eqref{eq:ori} with $f(X) = {\rm tr}(X\tp A X)$ and $v=e$. The precision of hashing is inherently influenced by various pre-processing, post-processing, and encoding procedures. Despite this consideration, our primary focus is on assessing the methods' proficiency in solving binary orthogonal problems. To do this, we assess the methods' performance in solving \eqref{eq:unsup_dh} with randomly generated data. Specifically, the Laplacian matrix is generated by $A = I_n - Z\Lambda Z\tp$ with $Z=$\texttt{randn(n,m)}, \texttt{m=500}, and $\Lambda={\rm diag}(Z\tp e)\in\R^{m\times m}$. We employ following evaluation metrics:  CPU time and violation of the balance/orthogonal constraint. The objective function value is not considered as an evaluation metric is because objective function value comparison is meaningless when the feasibility is not satisfied. The degrees of constraint violations are measured by $\|B\tp B-nI_r\|_F$ and $\|B\tp e\|_2$. Our evaluation involves experiments with varying sizes of the variable, which are shown in Table \ref{tab:exp_r} and \ref{tab:exp_n}.

In the set of experiments shown in Table \ref{tab:exp_r}, we evaluate the methods across different numbers of rows, while maintaining the number of columns of the variable at either 16 or 32. In the set of experiments shown in Table \ref{tab:exp_n}, we assess the methods across different numbers of columns, while keeping the number of rows of the variable at either 1,000 or 5,000.
The results indicate that the MGBO method consistently yields solutions with minimal violations of both the orthogonal and balanced constraints, outperforming other methods in most cases. Regarding  computational cost measured by CPU time, the MGBO method exhibits a slightly higher computational expense compared to the relaxation method, but it outperforms other methods in this aspect. 

\subsection{Results for the Supervised Hashing}\label{sec:sdh}

In the experiment on supervised discrete hashing, we test the following widely-used supervised hashing model\cite{shen2015supervised, shen2016fast, gui2017fast,gui2016supervised}:
\begin{align}\label{eq:sup_dh}
 &\underset{B\in\{-1,1\}^{n\times r}, W\in \R^{r \times m}}{\min} \widetilde{f}(X) := \|Y - BW\|_F^2 + \frac{\delta}{2}\|W\|_F^2 , \nonumber\\
 &\qquad\ {\rm s.t.~}\ B\tp B = nI_r,\quad B\tp e = 0,
\end{align}
where $\delta$ is a regularization parameter, $Y\in \R^{n\times m}$ is the label matrix with $m$ representing the number of label classes, and $W$ is the projection matrix that is learned jointly with $X$. The label matrix is sparse, with $Y_{ij}=1$ indicating that the $i$th sample belongs to the $j$th class, and $0$ otherwise. The projection matrix $W$, which is also present in the formulation, is learned simultaneously with $B$. We follow the implementation in \cite{gui2017fast, shen2015supervised, xiong2021generalized} to calculate the projection matrix $W$, which can be solved with a closed-form solution:
\[
W=(B\tp B + \delta I)^{-1}B\tp Y.
\]

We adopt  the implementation in \cite{gui2017fast,shen2016fast, shen2015supervised, xiong2021generalized} and use the linear hash function for encoding queries into binary code. The experiments are conducted on the MNIST digit dataset and the CIFAR-10 dataset, which are described as follows.
\begin{itemize}
 \item[i.] The \textbf{MNIST}\footnote{MINIST dataset https://pjreddie.com/projects/mnist-in-csv.} digit dataset consists of 70,000 samples of handwritten 0-9 digits. These images are all $28\times 28$ resolutions, which result in 784-dimensional feature vectors for representation. For the MNIST dataset, we split it into two subsets: a training set containing 69, 000 samples and a testing set set of 1,000 samples.
 \item[ii.] The \textbf{CIFAR-10}\footnote{CIFAR-10 dataset  https://www.cs.toronto.edu/~kriz/ cifar.html.} database is a collection of $32\times 32$ colour images consisting of five training batches and in 10 classes, with 6,000 images per class. There are 50,000 training images and 10,000 test images. The dataset is divided into five training batches and one test batch, each with 10,000 images. The training batches contain exactly 5,000 images from each class.
\end{itemize}

\begin{figure*}[h] 
\centering
\includegraphics[width=1.1\textwidth]{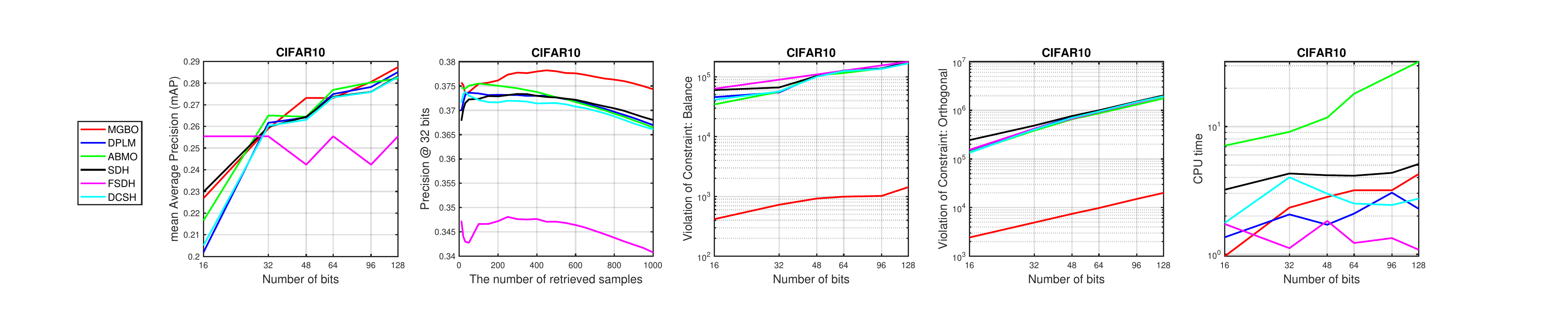}
\caption{\small Supervised hashing: CIFAR10 Dataset}
\label{fig:cifar}
\end{figure*}

\begin{figure*}[h] 
\includegraphics[width=1.1 \textwidth]{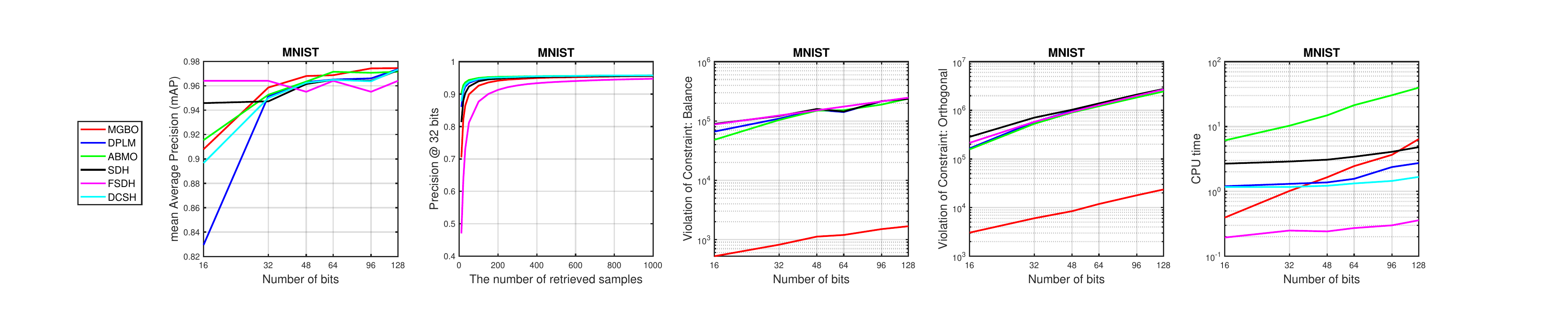}
\centering
\caption{\small Supervised hashing: MNIST Dataset}
\label{fig:mnist}
\end{figure*}

To evaluate performance, we compare the MGBO method with other hash solvers, including DPLM\cite{shen2016fast}, ABMO\cite{xiong2021generalized}, SDH\cite{gui2017fast}, FSDH\cite{gui2017fast}, and DCSH\cite{hoang2020unsupervised}.
We employ the Hamming ranking procedure for evaluations, which ranks database points based on their Hamming distances to the query vectors. Specifically, we measure their Hamming ranking precision using Mean Average Precision (MAP) and Precision. MAP calculates the average precision across all retrieved items for each bit, while Precision quantifies the proportion of relevant items among the top 1 to 1,000 retrieved items.
In addition, since binary hashing requires that the learned binary codes are uncorrelated and balanced in terms of the number of bits, we measure the performance of the solvers by assessing the feasibility violations of the returned solutions. Similar to the previous experiment,  $\|B\tp B-nI_r\|_F$ and $\|B\tp e\|_2$ are used to measure the degree of constraint violation, where the binary code $B$ is obtained by $B=\sgn(X^*)$, and $X^*$ is the solution returned by a solver.

We can observe from Figure \ref{fig:cifar} and \ref{fig:mnist} that the performance of the hash function associated to MGBO measured by MAP or precision is comparable to the performance of the hash functions associated to other solvers. Notably, MGBO outperform all other solvers in the cases of $r\in\{48, 128\}$ for the CIFAR10 dataset and $r\in\{48, 96, 128\}$ for the MNIST dataset.
MGBO also demonstrates superior performance in terms of feasibility violations compared to other methods (as shown in the third and the fourth sub-figures of \ref{fig:cifar} and \ref{fig:mnist}), which is consistent with the results of previous experiments. However, though MGBO provides a solution with minimal constraint violations, this does not guarantee its hash function is definitively better than others.

The FSDH method solves an optimization model without the binary constraint. As a result, it usually requires minimal CPU time for solving problem. In contrast,  aim to address both the binary constraint and equality constraints, incurring higher computational costs. Among these methods, ABMO records the highest computational cost. MGBO requires comparable computational costs with DPLM, SDH, and DCSH.

\section{Conclusions}\label{sec:conc}

 In this paper, we address the binary orthogonal optimization problem \eqref{eq:sh}. By investigating two penalized problems for its equivalent reformulation \eqref{eq:ori}, we have proposed two well-defined continuous optimization models \eqref{eq:penprob} and \eqref{eq:sep1}. It was shown that they serve as the global exact penalties under the assumption that the original model \eqref{eq:sh} is well-defined, i.e., has a nonempty feasible set. Based on the smooth penalty problem \eqref{eq:sep1}, we have presented a manifold gradient descent method with nonmonotone line-search rule, which is validated to have superior performance in terms of two evaluation indices and feasibility violation via supervised and unsupervised hashing tasks on MNIST and CIFAR-10 datasets. It is worth noting that the nonsmooth penalty problem \eqref{eq:penprob} can be solved using the manifold inexact augmented Lagrangian method \cite{deng2019inexact} which, after numerical testing, may yield solutions with a desirable feasibility violation after being projected onto the binary space, but will incur slightly higher computational costs and generate slightly inferior solutions compared to MGBO.

\section*{Acknowledgments}
 This work is supported by the National Natural Science Foundation of China under project No.11971177, Basic and Applied Basic Research Foundation of Guangdong Province under project No.2022A1515110959, and Guangzhou Municipal science and Technology Project under project No.202201010566.

\appendix

\section{Proofs}
In this section, we present the proofs of in this paper. 

\subsection{Proof of Lemma \ref{lemma3:bnr}}
\begin{proof}
We only consider the case when $v\neq 0$, because the equivalence $\bnr\equiv\snr$ holds when $v=0$, and the geometric properties of Stiefel manifold has been well-studied in many monographs such as \cite{boumal2020introduction, Absil08}. Let $F\!:\Rnr \rightarrow \mathbb{S}^r \times \R^r$ be a mapping defined by
 \[
 G(X):=\left(\begin{matrix}
 	X\tp X- I_r\\
 	X^{\top}v
 \end{matrix}\right)\quad{\rm for}\ X\in\mathbb{R}^{n\times r}.
 \]
 Obviously, $G$ is infinitely differentiable, and for each $Y\in\mathbb{R}^{n\times r}$, $G(Y)=0$ iff $Y\in\bnr$. By \cite[Definition~3.6]{boumal2020introduction}, it suffices to argue that for each $X\in\bnr$, the differential operator $DG(X)\!:\mathbb{R}^{n\times r}\to \mathbb{S}^r \times \R^r$ is surjective or equivalently its adjoint $[DG(X)]^*\!:\mathbb{S}^r \times \R^r\to\mathbb{R}^{n\times r}$ is injective. Fix any $X\in\bnr$. Pick any $(H,\xi)\in\mathbb{S}^r\times\mathbb{R}^r$ such that
 $[DG(X)]^*(H,\xi)=0$, i.e., $2XH+v\xi^{\top}=0$. Along with $X^{\top}X=I_{r}$, we have
 $2H+X^{\top}v\xi^{\top}=0$. Note that $X^{\top}v=0$. Then, $H=0$ and $v\xi^{\top}=0$.
 From $v\xi^{\top}=0$, we have $v\|\xi\|^2=0$, which along with $v\in\mathbb{R}^{n}\backslash\{0\}$ implies that $\xi=0$. Thus, we show that $[DG(X)]^*\!:\mathbb{S}^r \times \R^r\to\mathbb{R}^{n\times r}$ is injective, so $\bnr$ is a smooth embedded manifold of codimension $k=\frac{1}{2}r(r+1)+r$.
 The compactness of $\bnr$ is immediate by that of ${\rm St}(n,r)$.
\end{proof}

\subsection{Proof of Lemma \ref{lemma2-tangent}}
\begin{proof}
  Fix any $Z\in\mathbb{R}^{n\times r}$. Note that ${\rm proj}_{\mathbb{L}_1\cap\mathbb{L}_2}(Z)$ is the unique optimal solution $Y^*$ of the following strongly convex minimization problem
  \[
    \min_{Y\in\mathbb{R}^{n\times r}}\Big\{\frac{1}{2}\|Y-Z\|_F^2\ \ {\rm s.t.}\ \ \mathcal{A}_1(Y)=0,\,\mathcal{A}_2(Y)=0\Big\}.
  \]
  From the optimality of $Y^*$ to this minimization problem, there exist $G\in\mathbb{X}$ and $H\in\mathbb{Y}$ such that
  \begin{subequations}
   \begin{equation}\label{temp-equa0}	
  	Y^*-Z+\mathcal{A}_1^*(G)+\mathcal{A}_2^*(H)=0,
  \end{equation}
   \begin{equation}\label{temp-equa1}	
   -\mathcal{A}_1(Z)+\mathcal{A}_1\mathcal{A}_1^*(G)+\mathcal{A}_1\mathcal{A}_2^*(H)=0,
  \end{equation}
  \begin{equation}\label{temp-equa2}	
   -\mathcal{A}_2(Z)+\mathcal{A}_2\mathcal{A}_1^*(G)+\mathcal{A}_2\mathcal{A}_2^*(H)=0.
   \end{equation}
  \end{subequations}
  Combining \eqref{temp-equa1} with $\mathcal{A}_1\circ\mathcal{A}_2^*=0$ yields that
  $G=\mathcal{A}_1(\mathcal{A}_1\mathcal{A}_1^*)^{\dagger}(Z)$, and combining \eqref{temp-equa2} with $\mathcal{A}_2\circ\mathcal{A}_1^*=0$ yields that
  $H=\mathcal{A}_2(\mathcal{A}_2\mathcal{A}_2^*)^{\dagger}(Z)$. Substituting the two equalities into \eqref{temp-equa0} leads to
\begin{align*}
  Y^*=&~{\rm proj}_{\mathbb{L}_1\cap\mathbb{L}_2}(Z) \\
  =&~(\mathcal{I}-\mathcal{A}_1^*\mathcal{A}_1(\mathcal{A}_1\mathcal{A}_1^*)^{\dagger}-\mathcal{A}_2^*\mathcal{A}_2(\mathcal{A}_2\mathcal{A}_2^*)^{\dagger})(Z).
  \end{align*}
  This, by the arbitrariness of $Z$ in $\mathbb{R}^{n\times r}$, implies the desired conclusion.
\end{proof} 	

\subsection{Proof of Lemma \ref{equiv-model}}\label{app:equiv}
\begin{proof}
  Write $\Omega_1\!:=\big\{B\in\{-1,1\}^{n\times r}\,|\,B^{\top}B=nI_{r}\big\}$ and $\Omega_2:=\big\{B\in\Rnr\,|\,B^{\top}B=nI_{r},\, B\in[-E,E]\big\}$. We first verify that
  $\Omega_1=\Omega_2$. To this end, it suffices to argue that $\Omega_2\subset\Omega_1$. Pick any $B\in \Omega_2$. Since $B\tp B = n I_r$, we have $\|B_{\cdot j}\|^2=n$ for each $j\in[r]$. Suppose that $B\notin \{-1,1\}^{n\times r}$. There exists $(\overline{i},\overline{j})\in[n]\times[r]$ such that $|B_{\overline{i}\overline{j}}|<1$. Together with $B\in [-E,E]$, we have $\|B_{\overline{j}}\|^2<n$, a contradiction to $\|B_{\cdot j}\|^2=n$ for each $j\in[r]$. Hence, $\Omega_1=\Omega_2$. This implies that problem \eqref{eq:sh} is equivalent to
  \begin{equation}\label{equiv:sh1}
   \underset{B\in[-E,E]}{\min}\Big\{\widetilde{f}(B)\ \ {\rm s.t.}\ \ B^{\top}B=nI_{r},\,B^{\top}v=0 \Big\}
  \end{equation}
  in the sense that they have the same global and local optimal solution sets. By setting $X=\frac{1}{\sqrt{n}}B$, model \eqref{equiv:sh1} becomes \eqref{eq:ori}. Observe that each feasible point of \eqref{eq:sh} is a local optimal solution due to the discreteness of the set $\Omega_1$, so does each feasible point of \eqref{equiv:sh1}. The conclusion then follows.
 \end{proof}
 \subsection{Proof of Lemma \ref{lemma-boxStnr}}

  \begin{proof}
  From the proof of Lemma \ref{equiv-model}, it follows that $|\overline{X}_{ij}|=\frac{1}{\sqrt{n}}$ for each $(i,j)\in[n]\times[r]$. Let $\delta=\frac{1}{2\sqrt{n}}$. Pick any $X\!\in\mathbb{B}(\overline{X},\delta)\cap\Lambda$. Then, for each $(i,j)\in[n]\times[r]$, we have $|X_{ij}-\overline{X}_{ij}|\le\delta$, so that ${\rm sign}(X_{ij})={\rm sign}(\overline{X}_{ij})$. Hence,
  \begin{align*}
   &{\rm dist}\big(X,{\rm St}(n,r)\cap\Lambda\big)\le\|X-\overline{X}\|_F\le\sum_{j=1}^r\sum_{i=1}^n|X_{ij}-\overline{X}_{ij}|\nonumber\\
   &\le\sqrt{n}\sum_{j=1}^r\sum_{i=1}^n|X_{ij}-\overline{X}_{ij}|\big(2|\overline{X}_{ij}|-|X_{ij}-\overline{X}_{ij}|\big)\nonumber\\
   &=\sqrt{n}\sum_{j=1}^r\sum_{i=1}^n\Big[2|\overline{X}_{ij}||X_{ij}-\overline{X}_{ij}|-|X_{ij}-\overline{X}_{ij}|^2\Big]\nonumber\\
   &=\sqrt{n}\sum_{j=1}^r\sum_{i=1}^n\Big[-2\langle\overline{X}_{ij},X_{ij}-\overline{X}_{ij}\rangle-|X_{ij}-\overline{X}_{ij}|^2\Big]\nonumber\\
   &=\sqrt{n}\sum_{j=1}^r\sum_{i=1}^n\Big[\frac{1}{n}-|X_{ij}|^2\Big]=\sqrt{n}\sum_{j=1}^r\big(1-\|X_{\cdot j}\|^2\big)\nonumber\\
   &\le\sqrt{nr}\sqrt{\sum_{j=1}^r(1-\|X_{\cdot j}\|^2)^2}\le\sqrt{nr}\|X\tp X-I_r\|_F\nonumber\\
   &=\sqrt{nr}\|X\tp X-X\tp\overline{X}+X\tp\overline{X}-\overline{X}^{\top}\overline{X}\|_F\\
   &\le \sqrt{nr}(\|X\|+\|\overline{X}\|)\|X-\overline{X}\|_F\\
   &\le\sqrt{nr}(1+\delta){\rm dist}(X,\snr)\le2\sqrt{nr}{\rm dist}(X,\snr),
  \end{align*}
  where the third inequality is due to $2|\overline{X}_{ij}|-|X_{ij}\!-\!\overline{X}_{ij}|\ge\frac{1}{\sqrt{n}}$ for all $(i,j)\in[n]\times[r]$, the second equality is using ${\rm sign}(X_{ij})={\rm sign}(\overline{X}_{ij})$ and $|X_{ij}|\le\frac{1}{\sqrt{n}}=|\overline{X}_{ij}|$ for each $(i,j)\in[n]\times[r]$. The desired result then follows.
\end{proof}

\subsection{Proof of Proposition \ref{prop-BnrOmega}}
\begin{proof}
  By Lemma \ref{lemma-boxStnr}, there exists a constant $\delta'>0$ such that ${\rm dist}(Z,{\rm St}(n,r)\cap\Lambda)\le 2\sqrt{nr}{\rm dist}(Z,\snr)$ for all $Z\in\mathbb{B}(\overline{X},\delta')\cap\Lambda$. Fix any $X\!\in\mathbb{B}(\overline{X},{\delta'}/{2})$. Let $X_{\Lambda}\!:=\!{\rm proj}_{\Lambda}(X)$. 	Clearly, $\|X_{\Lambda}-\overline{X}\|_F\le 2\|X-\overline{X}\|_F\le\delta'$. Then, ${\rm dist}(X_{\Lambda},{\rm St}(n,r)\cap\Lambda)\le 2\sqrt{nr}{\rm dist}(X_{\Lambda},\snr)$. Notice that $\bnr\cap\Lambda\subset{\rm St}(n,r)\cap\Lambda$ are two discrete sets, where the discreteness of ${\rm St}(n,r)\cap\Lambda$ is implied by the proof of Lemma \ref{equiv-model}. There exists $\delta_1\in(0,\delta']$ such that for all $X\in\mathbb{B}(\overline{X},\delta_1)$,
  \[
   {\rm dist}(X,\bnr\cap\Lambda)=\|X\!-\!\overline{X}\|_F={\rm dist}(X,{\rm St}(n,r)\cap\Lambda).
  \]
 Consequently, for all $X\in\mathbb{B}(\overline{X},\delta_1)$, it holds that
 \begin{align*}
  &{\rm dist}\big(X,\bnr\cap\Lambda\big)={\rm dist}(X,{\rm St}(n,r)\cap\Lambda) \\
  &\le\|X-X_{\Lambda}\|_F+{\rm dist}\big(X_{\Lambda},{\rm St}(n,r)\cap\Lambda\big)\\
  &\le\|X-X_{\Lambda}\|_F+2\sqrt{nr}{\rm dist}(X_{\Lambda},{\rm St}(n,r))\\
  &\le\|X-X_{\Lambda}\|_F+2\sqrt{nr}\big[\|X-X_{\Lambda}\|_F+{\rm dist}(X,{\rm St}(n,r))\big]\\
  &\le (2\sqrt{nr}+1)\big[{\rm dist}(X,\Lambda)+{\rm dist}(X,\bnr)\big].
 \end{align*}
 This shows that the desired conclusion holds.
\end{proof}

\subsection{Proof of Corollary \ref{coro:gerror}}
\begin{proof}
 Fix any $Z\in\bnr$. By Proposition \ref{prop-BnrOmega}, for each $X\in\Gamma:=\bnr\cap\Lambda$, there exists $\delta_{X}$ such that
 \[
  {\rm dist}(Y,\Gamma)\le(2\sqrt{nr}\!+\!1)[{\rm dist}(Y,\Lambda)+{\rm dist}(Y,\bnr)].
 \]
 for all $Y\!\in\mathbb{B}(X,\delta_{X})$. Because $\bigcup_{X\in\Gamma}\mathbb{B}^{\circ}(X,\delta_{X})$ is an open covering of the compact set $\Gamma$, from Heine-Borel covering theorem, there exist $X^1,X^2,\ldots,X^p\in\Gamma$ such that $\Gamma\subset\bigcup_{i=1}^p\mathbb{B}^{\circ}(X^i,\delta_{X^i}):=D$. When $Z\in D$, from the last inequality, ${\rm dist}(Z,\Gamma)\le(2\sqrt{nr}\!+\!1){\rm dist}(Z,\Lambda)$, so we only need to consider that $Z\notin D$. Let $\overline{D}={\rm cl}[\bnr\backslash D]$. There must exist a constant $\widetilde{\kappa}>0$ such that $\min_{Y\in\overline{D}}{\rm dist}(Y,\Lambda)\ge\widetilde{\kappa}$. If not, there exists a sequence 	$\{Y^k\}_{k\in\mathbb{N}}\subset\overline{D}$ such that ${\rm dist}(Y^k,\Lambda)\le 1/k$, which by the compactness of the set $\overline{D}$ and the continuity of the distance function means that there is a cluster point, say $\overline{Y}\in\overline{D}$, of $\{Y^k\}_{k\in\mathbb{N}}$ such that $\overline{Y}\in\Lambda$. 	Then $\overline{Y}\in\Gamma\subset D$, a contradiction to the fact that $\overline{Y}\in\overline{D}$. In addition, because both $\overline{D}$ and $\bnr$ are compact, there exists a constant $c_0>0$ such that for all $Y\in\overline{D}$, ${\rm dist}(Y,\bnr)\le c_0$. Together with $\min_{Y\in\overline{D}}{\rm dist}(Y,\Lambda)\ge\widetilde{\kappa}$ and $Z\in\overline{D}$, it holds that ${\rm dist}(Z,\Gamma)\le (c_0/\widetilde{\kappa}){\rm dist}(Z,\Lambda)$. Thus, ${\rm dist}(Z,\Gamma)\le\kappa{\rm dist}(Z,\Lambda)$ with $\kappa=\max\{c_0/\widetilde{\kappa},2\sqrt{nr}\!+\!1\}$.
\end{proof}

\subsection{Proof of Theorem \ref{theorem1-epenalty}}
\begin{proof}
 According to Corollary \ref{coro:gerror}, ${\rm dist}\big(Z,\bnr\cap\Lambda\big)\le\kappa{\rm dist}(Z,\Lambda)$ for all $Z\in\bnr$. Fix any $X\in\bnr$. Then, there exists $\overline{X}\!\in\bnr\cap\Lambda$ such that
  \begin{equation}\label{eq:thm1_1}
	\|X\!-\!\overline{X}\|_F
	={\rm dist}(X,\bnr\cap\Lambda)
   \le \kappa{\rm dist}(X,\Lambda).
  \end{equation}
  Since $f$ is Lipschitz continuous on $\bnr$ with modulus $L_{\!f}$,
  \begin{equation}\label{eq:thm1_2}
	f(X)-f^*\ge f(X)-f(\overline{X})\ge-L_{\!f}\|X-\overline{X}\|_F.
  \end{equation}
 The first inequality is by the fact that $\overline{X}\!\in\bnr\cap\Lambda$, and the second inequality is by the Lipschitz continuity of $f$. By \eqref{eq:thm1_1} and \eqref{eq:thm1_2}, It follows that
  \begin{align}\label{eq:thm1_3}
    f(X)-f^*&\ge -\kappa L_{\!f}{\rm dist}(X,\Lambda)\ge- L_{\!f}\kappa{\rm dist}_1(X,\Lambda)\nonumber\\
    & =-\kappa L_{\!f}\,h_{\frac{1}{\sqrt{n}}}(X).
  \end{align}
  The second inequality is by the relation between dist$_1$ and dist. Thus, the first part of the conclusions follows. By \eqref{eq:thm1_3} and the compactness of $\bnr$, applying the result in  \cite[Proposition 2.1(b)]{liu2018equivalent}, we can conclude that problem \eqref{eq:penprob} associated to each $\rho\ge\kappa'L_f$ has the same global optimal solution set as problem \eqref{eq:ori} does.
\end{proof}

\subsection{Proof of Theorem \ref{thm:env_ep}}
\begin{proof}
 Let $X^*$ be a local optimal solution of \eqref{eq:ori}. There is $0<\varepsilon'<\frac{2}{\sqrt{n}}$ such that for all $X\in\B(X^*,\varepsilon')\cap(\bnr\cap\Lambda)$, it holds that
 \begin{equation}\label{eq:epthm1}
  f(X) \ge f(X^*).
 \end{equation}
 Let $\varepsilon=\frac{1}{2}\min\{\varepsilon',\delta',\gamma,\delta\}$, where $\delta$ is the same as the one in Proposition \ref{prop-BnrOmega}. Consider any $X\in\B(X^*,\varepsilon)\cap\bnr$.
 Pick any $\overline{X}\in{\rm proj}_{\bnr\cap\Lambda}(X)$, we have $\|\overline{X}-X\|_F\le \|X-X^*\|_F$. Then,
 \[
   \|\overline{X}-X^*\|_F\le\|\overline{X}-X\|_F+\|X-X^*\|_F\le 2\|X^*-X\|_F\le\varepsilon,
 \]
 which along with \eqref{eq:epthm1} implies that $f(\overline{X})\ge f(X^*)$.
 Together with the given condition in \eqref{eq:epthm0}, it then follows that
 \begin{equation}\label{ineq1-fstar}
 	f(X) - f(X^*)\ge f(X) - f(\overline{X})\ge - L'\|X - \overline{X}\|_F^2.
 \end{equation}
 In addition, from Proposition \ref{prop-BnrOmega}, it follows that
 \begin{equation}\label{ineq2-fstar}
   \|X-\overline{X}\|_F^2\le(2\sqrt{nr}\!+\!1)^2{\rm dist}^2(X,\Lambda).
 \end{equation}
 Notice that $\|X-\overline{X}\|_F\le\|X-X^*\|_F\le\varepsilon<\frac{1}{\sqrt{n}}$, while $|\overline{X}_{ij}|=\frac{1}{\sqrt{n}}$ for each $(i,j)\in[n]\times[r]$. It is immediate to deduce that ${\rm sgn}(X_{ij})={\rm sgn}(\overline{X}_{ij})$ for each $(i,j)\in[n]\times[r]$.
 Also, from $\|X-\overline{X}\|_F\le\varepsilon$, we have
 $\frac{1}{\sqrt{n}}-\varepsilon\le X_{ij}\le\frac{1}{\sqrt{n}}+\varepsilon$ for all $(i,j)\in[n]\times[r]$, which by $\varepsilon<\gamma$ implies that the index set $J_0:=\{(i,j)\in[n]\times[r]\ |\ |X_{ij}|>\frac{1}{\sqrt{n}}+\gamma\}$ is empty.
 Define $J_1:=\{(i,j)\in[n]\times[r]\ |\ |X_{ij}|<\frac{1}{\sqrt{n}}\}$ and
 $J_2:=\{(i,j)\in[n]\times[r]\ |\ \frac{1}{\sqrt{n}}\le|X_{ij}|\le\frac{1}{\sqrt{n}}+\gamma\}$.
 Together with Lemma \ref{lem:prox} with $c=\frac{1}{\sqrt{n}}$, it follows that
 \[
  {\rm env}_{\gamma h_{\frac{1}{\sqrt{n}}}}(X)
  =\frac{1}{2\gamma}\sum_{(i,j)\in J_2}|\overline{X}_{ij}-X_{ij}|^2.
 \]
 On the other hand, since $J_0=\emptyset$, it holds that
 \begin{align*}
   {\rm dist}^2(X,\Lambda)&=\sum_{(i,j)\in J_0\cup J_1\cup J_2}|[{\rm proj}_{\Lambda}(X)]_{ij}-X_{ij}|^2\\
   &=\sum_{(i,j)\in J_1}|X_{ij}-X_{ij}|^2+\sum_{(i,j)\in J_2}|\frac{1}{\sqrt{n}}{\rm sgn}(\overline{X}_{ij})-X_{ij}|^2 \\
   &=\sum_{(i,j)\in J_2}|\overline{X}_{ij}-X_{ij}|^2
 \end{align*}
 Henceforth ${\rm env}_{\gamma h_{\frac{1}{\sqrt{n}}}}(X)= \frac{1}{2\gamma}{\rm dist}^2(X,\Lambda)$. Together with inequality \eqref{ineq2-fstar},
 we obtain that
 \[
   {\rm env}_{\gamma h_{\frac{1}{\sqrt{n}}}}(X)\ge\frac{1}{2\gamma(2\sqrt{nr}+1)^2}\|X-\overline{X}\|_F^2.
 \]
 By combining above inequality with \eqref{ineq1-fstar}, it follows that
 \begin{equation}\label{eq:thm2_1}
   f(X)-f(X^*)+2L'\gamma(2\sqrt{nr}\!+\!1)^2{\rm env}_{\gamma h_{\frac{1}{\sqrt{n}}}}(X)
   \ge 0.
 \end{equation}
 The first part of the conclusions. By \eqref{eq:thm2_1} and the compactness of $\bnr$, applying the result in  \cite[Proposition 2.1(b)]{liu2018equivalent}, we can conclude that problem \eqref{eq:penprob} associated to each $\rho\ge2L'\gamma(2\sqrt{nr}\!+\!1)^2$ has the same global optimal solution set as problem \eqref{eq:ori} does.
\end{proof}

\subsection{Proof of Theorem \ref{thm:itr}}

\begin{proof}
Since $f$ is continuously differentiable on the compact manifold $\bnr$ and $h_{\frac{1}{\sqrt{n}}}(X)$ is bounded below, we have $\Theta_{\!\rho,\gamma}$ is also bounded below. Then there exists $\Theta_{\rho,\gamma}^*>-\infty$ such that for all $X\in \bnr$, $\Theta_{\!\rho,\gamma}(X)\ge \Theta_{\rho,\gamma}^*$. Suppose that Algorithm \ref{ManPG1} does not terminate after $K>(m+1)(2(\Theta_{\rho,\gamma}(X^0) - \Theta_{\rho,\gamma}^*)/(\alpha\eta \bar{t}\epsilon^2)+1)$ iterations, i.e., $\|{\rm grad}\,\Theta_{\!\rho,\gamma}(X^k)\|_F>\epsilon$ $\forall k=0,1,\dots, K-1$. Let $t_k$ be the stepsize in the $k$-th iteration, then $t_k\ge \eta\bar{t}$. Since $V^k = -t_k\,{\rm grad}\,\Theta_{\!\rho,\gamma}(X^k)$ and $X^{k+1} = R_{X^k}(V^k)$, we have
\begin{align*}
	\Theta_{\rho,\gamma}(X^{\ell(K)}) &\le \Theta_{\rho,\gamma}(X^{\ell(\ell(K)-1)})-\frac{\alpha}{2t_{\ell(K)-1}}\|V^{\ell(K)-1}\|_F^2\\
	&\le\Theta_{\rho,\gamma}(X^{\ell(K-m-1)})-\frac{\alpha}{2t_{\ell(K)-1}}\|V^{\ell(K)-1}\|_F^2\\
	&\le\Theta_{\rho,\gamma}(X^{0})-\sum_{i=0}^{\lfloor\frac{K}{m+1}\rfloor-1}\frac{\alpha\|V^{\ell(K\!-\!i(m\!+\!1))\!-\!1}\|_F^2}{2t_{\ell(K\!-\!i(m\!+\!1))\!-\!1}}\\
	&\le\Theta_{\rho,\gamma}(X^{0})-\frac{\alpha\epsilon^2}{2}\sum_{i=0}^{\lfloor\frac{K}{m+1}\rfloor-1}t_{\ell(K\!-\!i(m\!+\!1))\!-\!1}\\
	&\le\Theta_{\rho,\gamma}(X^{0})-\frac{\alpha\epsilon^2\lfloor\frac{K}{m+1}\rfloor}{2}\eta\bar{t}.
\end{align*}
Thus, we have $\Theta_{\rho,\gamma}(X^0) - \Theta_{\rho,\gamma}^* \ge\frac{\alpha\epsilon^2\lfloor\frac{K}{m+1}\rfloor}{2}\eta\bar{t}$
which means,
\[
	K \le (m+1)(\frac{2(\Theta_{\rho,\gamma}(X^0) - \Theta_{\rho,\gamma}^*)}{\alpha\eta \bar{t}\epsilon^2}+1).
\]
This yields a contradiction to the fact that $K> (m+1)(2(\Theta_{\rho,\gamma}(X^0) - \Theta_{\rho,\gamma}^*)/(\alpha\eta \bar{t}\epsilon^2)+1)$. The conclusion then follows.
\end{proof}

\end{document}